\newif\ifconfver
\confverfalse      

\newif\ifcutshort      
\cutshorttrue

\newif\ifcutshortlvltwo  

\ifconfver
    \documentclass[10pt,twocolumn,twoside]{IEEEtran}
\else
    \documentclass[11pt,draftcls,onecolumn]{IEEEtran}
\fi

\usepackage[dvips]{epsfig}
\usepackage[dvips]{graphicx}
\usepackage{boxedminipage}
\usepackage{color}
\usepackage{amsfonts}
\usepackage{amssymb}
\usepackage{amsmath}
\usepackage{graphics,amsmath,amsfonts,amssymb,epsfig,latexsym,amsfonts,graphicx,amssymb}
\usepackage{setspace}
\usepackage{cite}
\usepackage{wrapfig}
\usepackage{subfigure, psfrag, framed}
\usepackage{algpseudocode}
\usepackage{algorithmicx}

\usepackage{stfloats}
\usepackage{url}
\usepackage{subfig}

\newcommand{\bfPhi}{{\mbox{\boldmath $\Phi$}}}

\newcommand{\st}{{\rm s.t.}}

\newcommand{\bx}{\mathbf{x}}

\newcommand{\hu}{\widehat{u}}

\newcommand{\cI}{\mathcal{I}}
\newcommand{\cX}{\mathcal{X}}

\newcommand{\hx}{\widehat{x}}

\newtheorem{thm}{Theorem}

\newtheorem{coro}{Corollary}

\newtheorem{exa}{Example}

\begin{document}
\title{A Unified Algorithmic Framework for Block-Structured Optimization Involving Big Data}
\author{Mingyi Hong$^*$ \thanks{M. Hong is with the Department of Industrial and Manufacturing Systems Engineering, Iowa State University, Ames, IA.},  Meisam Razaviyayn$^*$ \thanks{M. Razaviyayn is with the Department of Electrical Engineering, Stanford University, CA.},  Zhi-Quan Luo, \thanks{Z.-Q. Luo is with the Department of Electrical and Computer Engineering, University of Minnesota, MN.} and Jong-Shi Pang\thanks{J.-S. Pang is with the Daniel J. Epstein Department of Industrial and Systems Engineering,
University of Southern California, CA.}
\thanks{* The first two authors contributed equally to this work. M. Hong is supported by NSF, Grant No. CCF-1526078. Z.-Q. Luo is supported by NSF, Grant No. CCF-1526434}} \maketitle
\begin{abstract}
This article presents a powerful algorithmic framework for big data optimization, called the Block Successive Upper bound Minimization (BSUM). The BSUM includes as special cases many well-known methods for analyzing massive data sets, such as the Block Coordinate Descent (BCD), the Convex-Concave Procedure (CCCP), the Block Coordinate Proximal Gradient (BCPG) method, the Nonnegative Matrix Factorization (NMF), the Expectation Maximization (EM) method and so on. In this article, various features and properties of the BSUM are discussed from the viewpoint of design flexibility, computational efficiency, parallel/distributed implementation and the required communication overhead. Illustrative examples from networking, signal processing and machine learning are presented to demonstrate the practical performance of the BSUM framework.

\end{abstract}

%

\section{Introduction}

\subsection{Overview of optimization for big data} \label{sub:bigdata}
With advances in sensor, communication and storage technologies, data acquisition is more ubiquitous than any time in the past.  This has made available big data sets in many areas of engineering, biological, social, and physical sciences. While the proper modeling and analysis of such data sets can yield valuable information for inference, estimation, tracking, learning and decision-making, their size and complexity present great challenges in algorithm design and implementation.

 Due to its central role in big data analytics, large-scale optimization has recently attracted significant attention not only from the optimization community, but also from the machine learning, statistics as well as the signal processing communities. For example, emerging problems in image processing, social network analysis and computational biology can easily exceed millions or billions of variables, and significant research is underway  to enable fast accurate solution of these problems {\color{black}\cite{Krishnamurthy14,bigdata14,Sra12,roberts2013streaming}} . Also, problems related to the design and provision of large scale smart infrastructures such as wireless communication networks require real-time efficient resource allocation decisions to ensure optimal network performance. Traditional general purpose optimization tools are inadequate for these problems due to the complexity of the model, the heterogeneity of the data, and most importantly the sheer data size {\cite{nestrov12,Cevher14, Aybat14, richtarik12,cui15}}. Modern large-scale optimization algorithms, especially those that are  capable of exploiting problem structures, dealing with distributed, time varying and incomplete data sets, utilizing massively parallel computing and storage infrastructures, have become the workhorse in the big data era. 

 To be efficient for big data applications, optimization algorithms must have certain properties:
 \begin{enumerate}
 \item Each of their computational steps must be simple and easy to perform;
 \item The intermediate results are easily stored;
 \item They can be implemented in a {distributed} and/or {parallel} manner so as to exploit the modern multi-core and cluster computing architectures;
 \item A high-quality solution can be found using a small number of iterations.
 \end{enumerate}
 These requirements preclude the use of high-order information about the problem (i.e., the Hessian matrix of the objective), which is usually too expensive to obtain, even for modest-sized problems. 

\subsection{The Block Coordinate Decent Method}\label{sub:bcd}
\begin{figure}
 \centering
 \vspace{-0.4cm}
\includegraphics[width=0.5\linewidth]{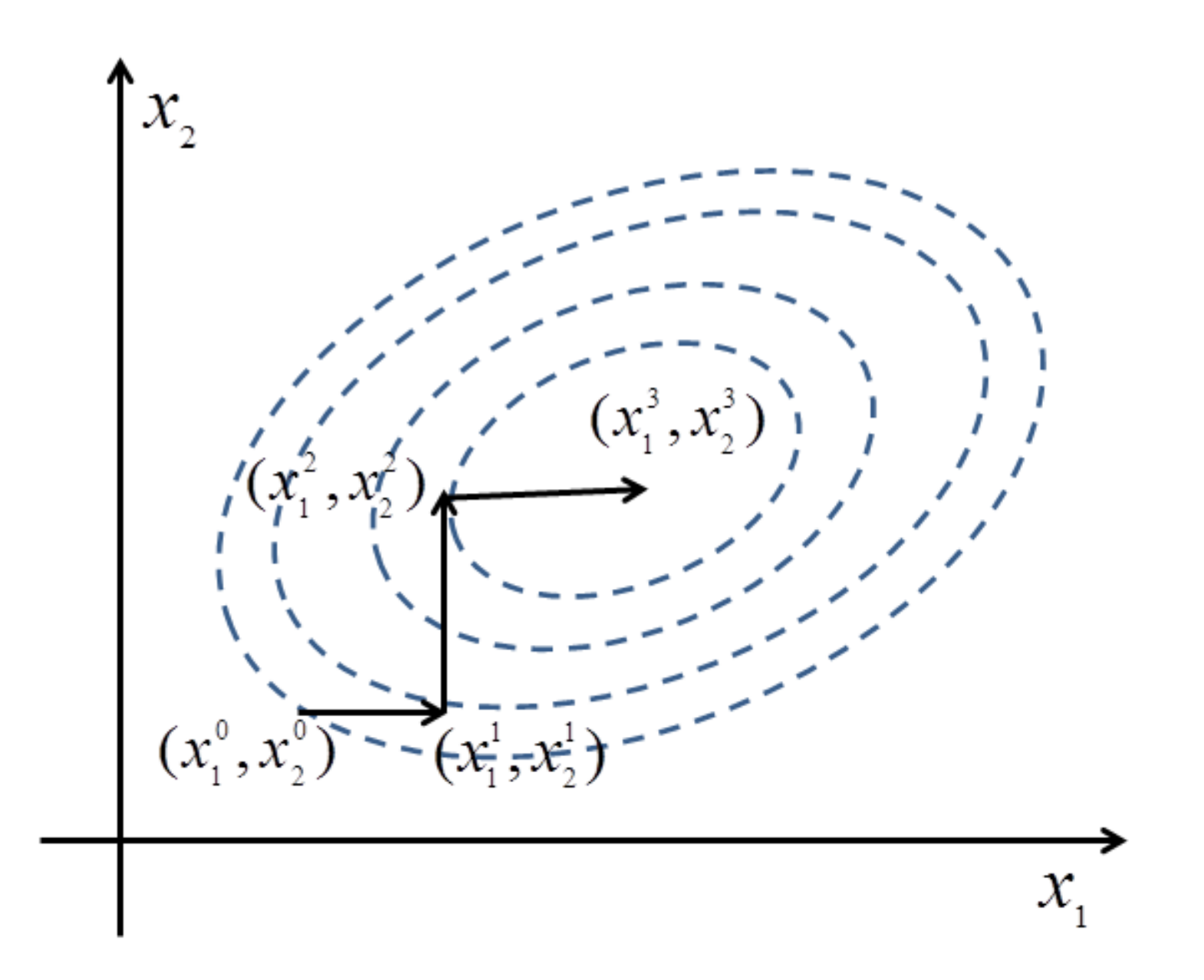}
\vspace*{-15pt}
\caption{\footnotesize Illustration of the BCD method for a two-dimensional problem. The dashed curves represent the contours of the objective function, the solid lines represent the progress of the algorithm.}\label{fig:BCD}
\end{figure}

 A very popular family of optimization algorithms that satisfies most of the aforementioned properties is the {\it block coordinate descent} (BCD) method, sometimes also known as the alternating minimization/maximization (AM) algorithm. The basic steps of the BCD are simple: {\bf i)} partition the entire optimization variables into small blocks, and {\bf ii)} optimize one block variable (or few blocks of variables) at each iteration, while holding  the remaining variables fixed. More specifically, consider the following block-structured optimization problem
 \begin{align}
\label{eq:OriginalProblem}
\displaystyle{ {\operatornamewithlimits{\mbox{minimize}}_{x}}} \quad  f(x_1,x_2,\ldots,x_n),\quad
\st \ x_i\in\mathcal{X}_i, \ i=1,2,...,n,
\end{align}
where $f(\cdot)$ is a continuous function (possibly nonconvex, nonsmooth), each $\mathcal{X}_i\subseteq \mathbb{R}^{m_i}$ is a closed
convex set,  and each $x_i$ is a block variable, $i=1,2,\ldots,n$. Define $x:=(x_1,\cdots, x_n)\in\mathbb{R}^m$, and let $\mathcal{X}:=\cX_1\times\cdots\times \cX_n\subseteq \mathbb{R}^{m}$. When applying the classical BCD method to solve \eqref{eq:OriginalProblem}, at every iteration $r$, a single block
of variables, say $i = (r \;{\rm
mod}\;n)+1$, is optimized by solving the following problem
\begin{equation}
\label{eq:GSperblock}
\begin{split}
x_i^r \in \arg\min_{x_i\in \mathcal{X}_i} \quad &f(x_i,x_{-i}^{r-1}),
\end{split}
\end{equation}
 where we have defined $x^{r-1}_{-i}:=(x^{r-1}_1,\cdots, x^{r-1}_{i-1},x^{r-1}_{i+1},\cdots, x^{r-1}_n)$; for the rest of variables $j\ne i$, let $x^r_j=x^{r-1}_j$. We refer the readers to Fig. \ref{fig:BCD} for a graphical illustration of the algorithm.

 The BCD algorithm is intuitively appealing and  very simple to implement, yet it is extremely powerful. It enjoys tremendous popularity in a wide range of applications from signal processing, communications, to machine learning.
 Representative examples include  image deblurring \cite{you1996regularization}, statistical learning \cite{Friedman10}, wireless communications \cite{shi11WMMSE_TSP}, and so on. In recent years, there is a surge of renewed interests to extend and generalize the BCD type of algorithms due to its applications in modern big data optimization problems. Theoretically, the BCD algorithm and its variants have been significantly generalized to accommodate various coordinate update rules \cite{Chen2012MBI, nestrov12, schmidt15icml}, and have been made suitable for  implementation on modern parallel processing architecture \cite{rick13parallel,scutari13flexible, liu14asynchronous,meisam14nips,Facchinei15}. It can handle a wide range of nonsmooth nonconvex cost functions \cite{razaviyayn2013stochastic, Razaviyayn12SUM, Mairal13}.  Practically, it has been used in emerging large-scale signal processing and machine learning applications, such as compressive sensing/sparse signal recovery \cite{Beck:2009:FIS:1658360.1658364, Figueiredo07}, matrix completion \cite{candes:2012:EMC:2184319.2184343, sunluo14}, tensor decomposition \cite{kim2013algorithms, TensorReviewKolda}, to name just a few. A recent survey of this algorithm can be found in \cite{wright15cd}.


\subsection{The Block Successive Upper Bound Minimization Method}\label{sub:bsum}
In this article, we introduce a unifying framework, called the Block Successive Upper bound Minimization (BSUM) method, which generalizes the BCD type of algorithms \cite{Razaviyayn12SUM}. Simply put, the BSUM framework includes algorithms that successively optimize certain upper-bounds or surrogate functions of the original objectives, possibly in a block by block manner. The BSUM framework significantly expands the application domain of the traditional BCD algorithms. For example, it covers many classical statistics and machine learning algorithms such as the Expectation Maximization (EM) method \cite{EMDempster}, the Concave-Convex Procedure (CCCP) \cite{CCCP} and the multiplicative Nonnegative Matrix Factorization (NMF) \cite{Lee00algorithmsfor}. It also includes as special cases many well-known signal processing algorithms such as the family of interference pricing algorithms \cite{Shi:2009, hong12survey} and the weighted minimum mean square error (WMMSE) algorithms \cite{shi11WMMSE_TSP, hong12HetNetFramework} for interference management in large-scale wireless systems.

The generality and flexibility of the BSUM offers an ideal platform to explore optimization algorithms for big data. Through the lens of the BSUM, one obtains not only a thorough understanding of the variety of algorithms and  applications that are being covered, but more importantly the principle for designing a good algorithm suitable for big data. {To this end, this article will first provide a concise overview of a few key theoretical characterizations of the algorithms that fall in the BSUM framework, BCD included.}  {Our emphasis will be given to providing intuitive understanding as to when and where the BSUM framework should (or should not) work, and how its performance can be characterized.}  {The second part of this article offers a detailed account of many existing large-scale optimization algorithms that fall in the BSUM framework, together with a few big data related applications that are of significant interests to the signal processing community. The last part of the article outlines some interesting extensions of the BSUM that further help expand its application domains. {Throughout this article, special emphasis will be given to computational issues arising from big data optimization problems such as algorithm design for parallel and distributed computation, algorithm implementation on multi-core computing platforms and distributed storage of the data. In particular, we will discuss how these issues can be addressed in the BSUM framework by providing theoretical insights and examples from real practical problems.

\section{BSUM and Its Theoretical Properties}
\subsection{The Main Idea}
We start by providing a high level description of the BSUM method and some of its theoretical properties. 
In practice one of the main problems of directly applying BCD to solve problem \eqref{eq:OriginalProblem} is that each of its subproblem \eqref{eq:GSperblock} is often difficult to solve exactly, especially when $f(x)$ is nonconvex. Moreover, even such exact minimization can be performed, the BCD may not necessarily converge. One of the key insights to be offered by the BSUM framework is that, for both practical and theoretical considerations, obtaining an {\it approximate} solution of \eqref{eq:GSperblock} is good enough to keep the algorithm going. To be more specific, let us introduce $u_i(x_i, z): \cX_i\to \mathbb{R}$ as an {\it approximation function} of $f(x_i, z_{-i})$ for each coordinate $i$ at a given feasible point $z\in \cX$.
Let us define a set $\cI^{r}$ (possibly with $|\cI^r|>1$) as the block-variable indices to be picked at iteration $r$. Then at each iteration $r$, the BSUM method performs the following simple update
\begin{align}\label{eq:BSUMperblock}
\left\{
\begin{array}{l}
x_i^{r}\in {\rm arg}\!\min_{x_i\in \mathcal{X}_i}u_i(x_i, x^{r-1}),\; \forall~i\in\cI^r\\
x_k^{r}=x^{r-1}_k,\ \forall~k\notin \cI^r
\end{array}
\right..
\end{align}

\begin{table}[htb]
\centering
\begin{tabular}{|p{4in}|}
\hline
\begin{itemize}
\item [1] \;Find a feasible point $x^0\in {\cX}$ and set $r = 0$
\item [2] \; \textbf{repeat}
\item [3] \quad {Pick index set $\cI^{r}$}
\item [4] \quad  {Let $x_i^{r}  \in \arg \min_{x_i\in \mathcal{X}_i} u_i(x_i,x^{r-1})$, $\forall~i\in\cI^r$}
\item [5] \quad Set $x_k^{r} = x_k^{r-1}, \quad \forall\; k \notin \cI^r$
\item [6] \quad $r = r+1$,
\item [7] \; \textbf{until} some convergence criterion is met
\end{itemize}
\\
\hline
\end{tabular}\vspace{1.2em}
\caption{Pseudo code of the BSUM algorithm}
\label{Table:BSUMAlgorithm} \vspace{-0.5cm}
\end{table}
The complete description of the BSUM is given in Table \ref{Table:BSUMAlgorithm}. We also refer the readers to Fig. \ref{fig:BSUM} for a graphical comparison of the iterates generated by BSUM and BCD for a two dimensional problem. It should be clear at this point that when $\cI^r = \{(r \;{\rm mod} \;n)+1\}$ and
no approximation is used (i.e., $u_i(x_i, z)=f(x_i, z_{-i})$ and at each iteration a single coordinate is selected), then the BSUM reduces to the classical cyclic BCD method. In Table \ref{table:coordinate} we also present several index selection rules that are covered by the BSUM framework. For simplicity of presentation, we will use the classical cyclic index selection rule where $\cI^r = \{(r \;{\rm mod} \;n)+1\}$ in the remainder of this article unless otherwise noted.

 \begin{table*}[h]
\centering
\begin{tabular}{|p{6.5in}|}
\hline\\

At each iteration $r$, define a set of auxiliary variables $\{\hx_i^{r}\}_{i=1}^{n}$ as:
\begin{align}
\hx^{r}_i\in \arg\min_{x_i\in \cX_i}\; u_i\left(x_i, x^{r-1}\right),  \; i=1,\cdots, n. \nonumber
\end{align}
Then we have the following commonly used coordinate selection rules.
\begin{itemize}
\item {\em Cyclic rule}: The coordinates are chosen cyclically, i.e., in the order of $1,2, \cdots, n, 1, 2, \cdots$.
\item {\em Essentially cyclic (E-C) rule}: There exists a given period $T\ge 1$ during which each block is updated at least once, i.e.,
$$\bigcup_{i=1}^{T}\cI^{r+i}=\{1,\cdots,n\}, \; \forall~r.$$

\item  {\it Gauss-Southwell (G-So) rule}: At each iteration $r$, $\cI^{r}$ contains a single index $i^*$ that satisfies:
$$i^*\in \left\{i \;\bigg{|}\; \|\hx^{r}_{i}-x^{r-1}_{i}\|\ge q \max_{j}\|\hx^{r}_j-x^{r-1}_j\|\right\}
$$ for some constant $q\in (0,\; 1]$.

\item {\it Maximum block improvement (MBI) rule}: At each iteration $r$, $\cI^{r}$ contains a single index $i^*$ that achieves the best objective: $i^*\in \arg\min_{i}f(\hx^{r}_i, x^{r-1}_{-i}).$

\item {\it Randomized rule}: {Let $p_{\min}\in (0,1)$ be a constant. At each iteration $r$, there exists a probability vector $p^{r} = (p_1^r, \ldots,p_n^r)\in\mathbb{R}^{n}$ satisfying $\sum_{i=1}^n p_i^r=1$ and $p_i^r>p_{\min}, \forall i$, with $\cI^r$ containing a single random index $i_r^*$ determined by
    $$\mbox{Pr}(i\in\cI^r \mid x^{r-1}, x^{r-2},\cdots, x^0)=p^r_i, \; \forall~i.$$}
\end{itemize}\\
\hline
\end{tabular}
\caption{{Commonly Used Coordinate Selection Rules.}}\label{table:coordinate}
\end{table*}

 \begin{figure*}[ht]
 \begin{minipage}[t]{0.48\linewidth}
\includegraphics[width=1\linewidth]{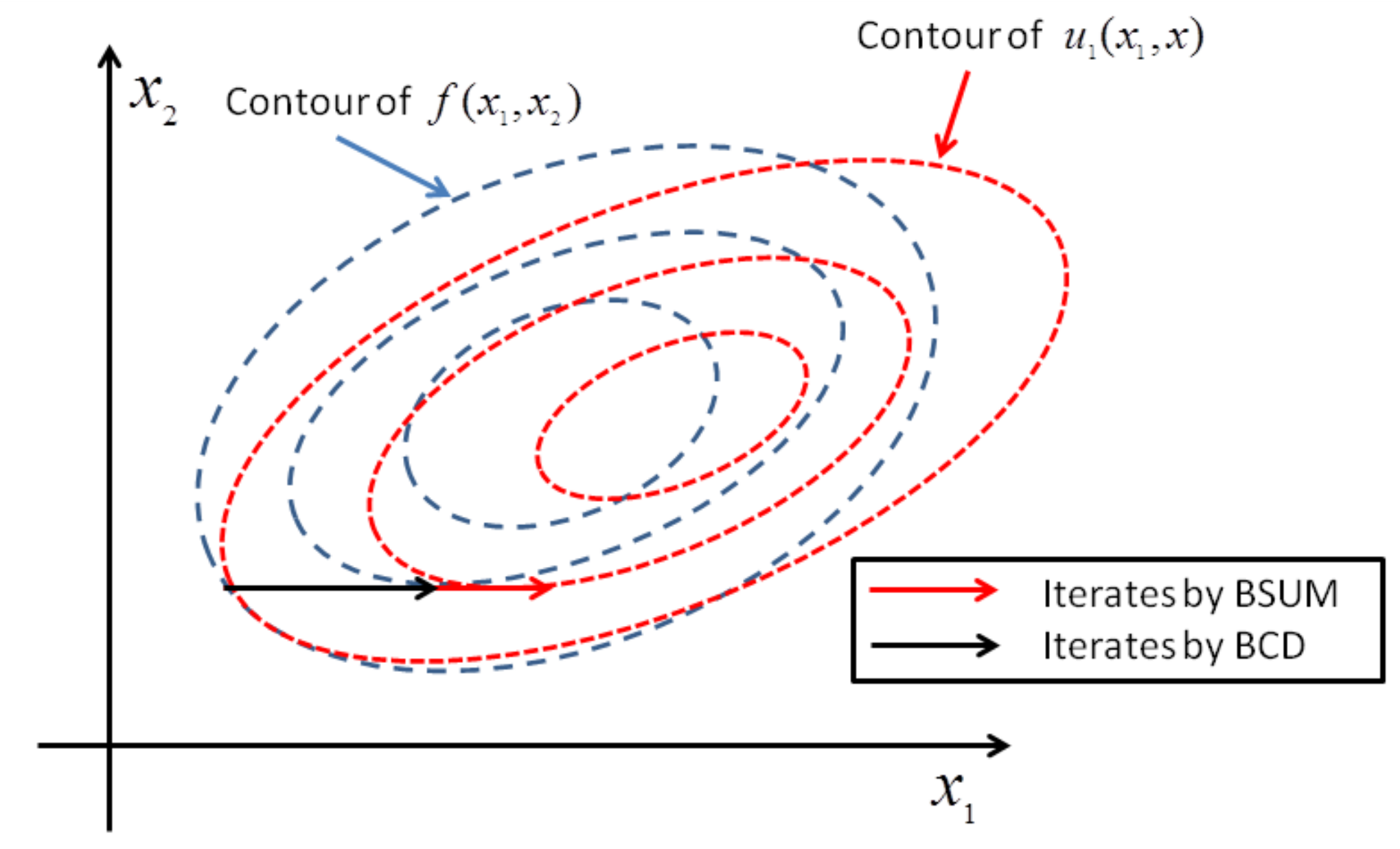}
\end{minipage}\hfill
 \begin{minipage}[t]{0.48\linewidth}
\includegraphics[width=1\linewidth]{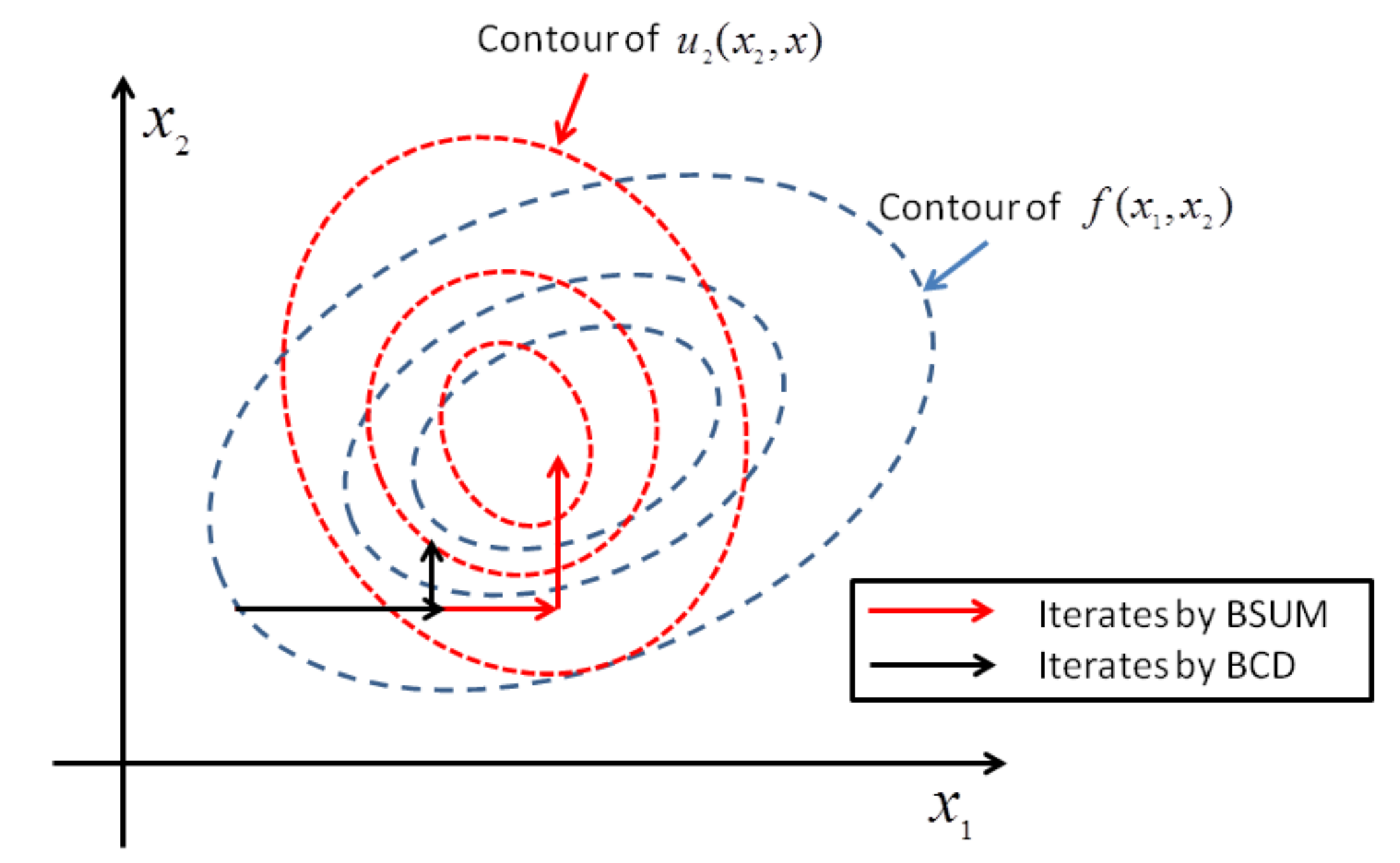}
\end{minipage}
\caption{\footnotesize Comparison of BSUM and BCD methods for solving a two-dimensional problem. Each time a single coordinate is picked for update. The blue dashed curves represent the contours of the objective function $f(x_1,x_2)$, the red dotted curves represent the contours of the upper bound functions $u_1(x_1,x)$ and $u_2(x_2,x)$; the black (resp. red) solid lines represent the progress of the BCD (resp. BSUM).}\label{fig:BSUM}
\end{figure*}

{Next we introduce the precise definition of the approximation function.} The main idea is that for each $i$, the approximation $u_i(x_i, x^r)$ should be an {\it upper bound} of the original objective function at the point of $x^r$ (hence the BSUM name of the framework). Please see Fig. \ref{fig:upperbound-step} for an illustration of the upper bound minimization process. Intuitively, picking an upper bound approximation function is reasonable because optimizing it at least should guarantee some descent of the original objective $f$; see Fig. \ref{fig:upperbound-step}(c).

{To be more precise, let us first define the directional derivative of a given function $f(x): \cX\to \mathbb{R}$ at a point $x\in \cX$ in direction $d$:
\[
f'(x;d) \triangleq \liminf_{\lambda \downarrow 0} \frac{f(x+\lambda d) - f(x)}{ \lambda}.
\]
Using this definition, we make the following assumptions on the $u_i$'s.

{\bf Assumption A.}
\begin{align}
& u_i(x_i, x) = f(x), \quad \forall\; x\in \mathcal{X}, \; \forall\; i \label{A1}\tag{A1} \\
& u_i(x_i, z) \geq f(x_i, z_{-i}),\quad\; \forall\; x_i \in \mathcal{X}_i, \; \forall\; z \in\mathcal{X}, \forall\; i\label{A2} \tag{A2}\\
&  u_i' (x_i,z;d_i)\bigg|_{x_i = z_i} =  f' (z;d), \quad \forall\; d = (0,\ldots,d_i,\ldots,0) \;\;\st \;\; z_i + d_i \in \mathcal{X}_i, \; \forall\; i \label{A3}\tag{A3}\\
& u_i(x_i, z) \; {\rm is} \; {\rm continuous} \; {\rm in} \; (x_i,z)
\label{A4}\tag{A4}, \quad \forall\; i.
\end{align}

 Intuitively, Assumptions \eqref{A1} and \eqref{A2} imply that the approximation function is a {\it global upper bound} of $f(x)$; while the assumption \eqref{A3} guarantees that the first order behaviors of the objective function and the approximation function are the same at the point of approximation (cf. Fig. \ref{fig:upperbound-step}). 
 In Table \ref{table:bounds}, we provide the readers with a few commonly used upper bounds that satisfy Assumption A;  see also \cite{Razaviyayn12SUM,Mairal13} for additional examples. More discussion will be given in subsequent sections on how these bounds are used in practice.

 \begin{table*}[h]
\centering
\begin{tabular}{|p{6.2in}|}
\hline
\begin{itemize}
\item {\em Proximal Upper Bound}: Given a constant $\gamma>0$, one can construct a bound by adding a quadratic penalization (i.e., the proximal term)
$$u_i(x_i,z) := f(x_i,z_{-i})+\frac{\gamma}{2}\|x_i-z_i\|^2.$$
\item {\em Quadratic Upper Bound}: Suppose $f(x)=g(x_1,\cdots, x_n)+h(x_1, \cdots, x_n)$, where $g$ is smooth with $\mathbf{H}_i$ as the Hessian matrix for the $i$th block. 
Then one can construct the following bound
$$u_i(x_i,z) := g(z_i,z_{-i})+h(x_i,z_{-i})+\langle \nabla_i g(z_i,z_{-i}), x_i-z_i\rangle+\frac{1}{2}(x_i-z_i)^T\bfPhi_i(x_i-z_i),$$
where both $\bfPhi_i$ and $\bfPhi_i- \mathbf{H}_i$ are positive semidefinite matrices.
\item {\em Linear Upper Bound}: Suppose $f$ is differentiable and concave, then one can construct
$$u_i(x_i,z) := f(z_i,z_{-i})+\langle \nabla_i f(z_i,z_{-i}), x_i-z_i\rangle.$$
\item {\em Jensen's Upper Bound}: Suppose $f(x):=f(a^T_1x_1,\cdots, a^T_n x_n)$ where $a_i\in\mathbb{R}^{m_i}$ is a coefficient vector, and $f$ is convex with respect to each $a^T_ix_i$. Let $w_i\in\mathbb{R}_{+}^{m_i}$ denote a weight vector with $\|w_i\|_1=1$. Then one can use Jensen's inequality and construct
$$u_i(x_i,z) := \sum_{j=1}^{m_i}w_i(j) f\left( \frac{a_i(j)}{w_i(j)}\left(x_i(j)-z_i(j)\right) +a^T_i z_i, z_{-i}\right)  $$
where $w_i(j)$ represents the $j$-th element in vector $w_i$.
\end{itemize}\\
\hline
\end{tabular}
\caption{{Commonly Used Upper Bounds Satisfying Assumption A.}}\label{table:bounds}
\end{table*}

It is worth mentioning that for a popular subclass of problem \eqref{eq:OriginalProblem}, Assumption A can be further simplified; see the following example \cite[Proposition 2]{Razaviyayn12SUM}.

\begin{exa}\label{ex:separability}
{\it Consider the following special form of problem \eqref{eq:OriginalProblem}:
 \begin{align}
\min f(x):=g(x_1,\cdots, x_n)+h(x_1,\cdots, x_n),\quad \st\;  \; x_i\in \cX_i, \ i=1,\cdots, n,\label{eq:problemBSUM1}
\end{align}
where $g: \cX\to \mathbb{R}$ is a smooth function and $h: \cX\to \mathbb{R}$ is a possibly nonsmooth function whose directional derivative exists at every point $x\in \cX$. Consider $u_i(x_i,z)=\hu_{i}(x_i,z)+h(x)$, where $\hu_{i}(x_i,z)$ approximates the smooth function $g$ in the objective. Then  assumption \eqref{A3} is implied by \eqref{A1} and \eqref{A2} and is therefore no longer needed.}
\end{exa}

\begin{table*}[h]
\centering
\begin{tabular}{|p{4.9in}|}
\hline
\begin{itemize}
\item {\em{Stationary solutions:}} The point~$x^*$
is a {\it stationary solution} of \eqref{eq:OriginalProblem} if $f'(x^*;d) \geq 0$ for all $d$
such that $x + d \in \cX$. Let $\cX^*$ denote the set of stationary solutions.
%

\item {\em{Coordinatewise minimum solutions}:} $\hat{x} \in \cX $ is coordinatewise minimum of problem \eqref{eq:OriginalProblem} with respect to the coordinates in $\mathbb{R}^{m_1},\mathbb{R}^{m_2},\ldots, \mathbb{R}^{m_n}$, if
     \[
     f(\hat{x} + d_k^0) \geq f(\hat{x}),\quad \forall\;  d_k \in \mathbb{R}^{m_k} \quad {\rm with} \quad \hat{x} +  d_k^0\in \cX,\; \;\forall k = 1,2,\ldots,n,
     \]
     where $d^0_k = (0,\ldots,d_k,\ldots, 0)$.
\end{itemize}\\
\hline
\end{tabular}
\caption{{Optimality Conditions.}}\label{table:def}
\end{table*}

\begin{figure}
\begin{center}
{ \subfigure[]
[]{\resizebox{.3\textwidth}{!}{\includegraphics{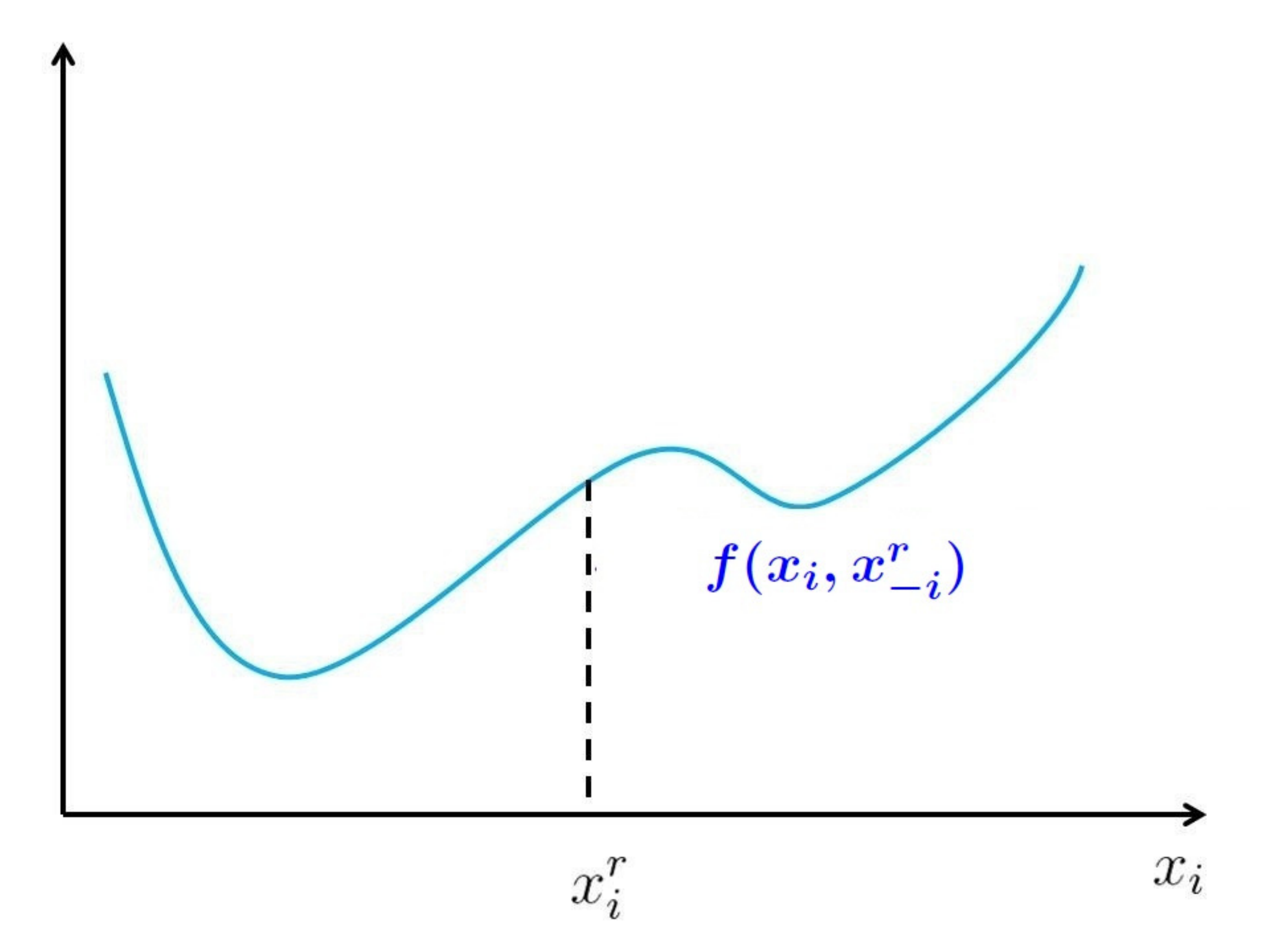}}}}
\hspace{1pc} { \subfigure[]
[]{\resizebox{.3\textwidth}{!}{\includegraphics{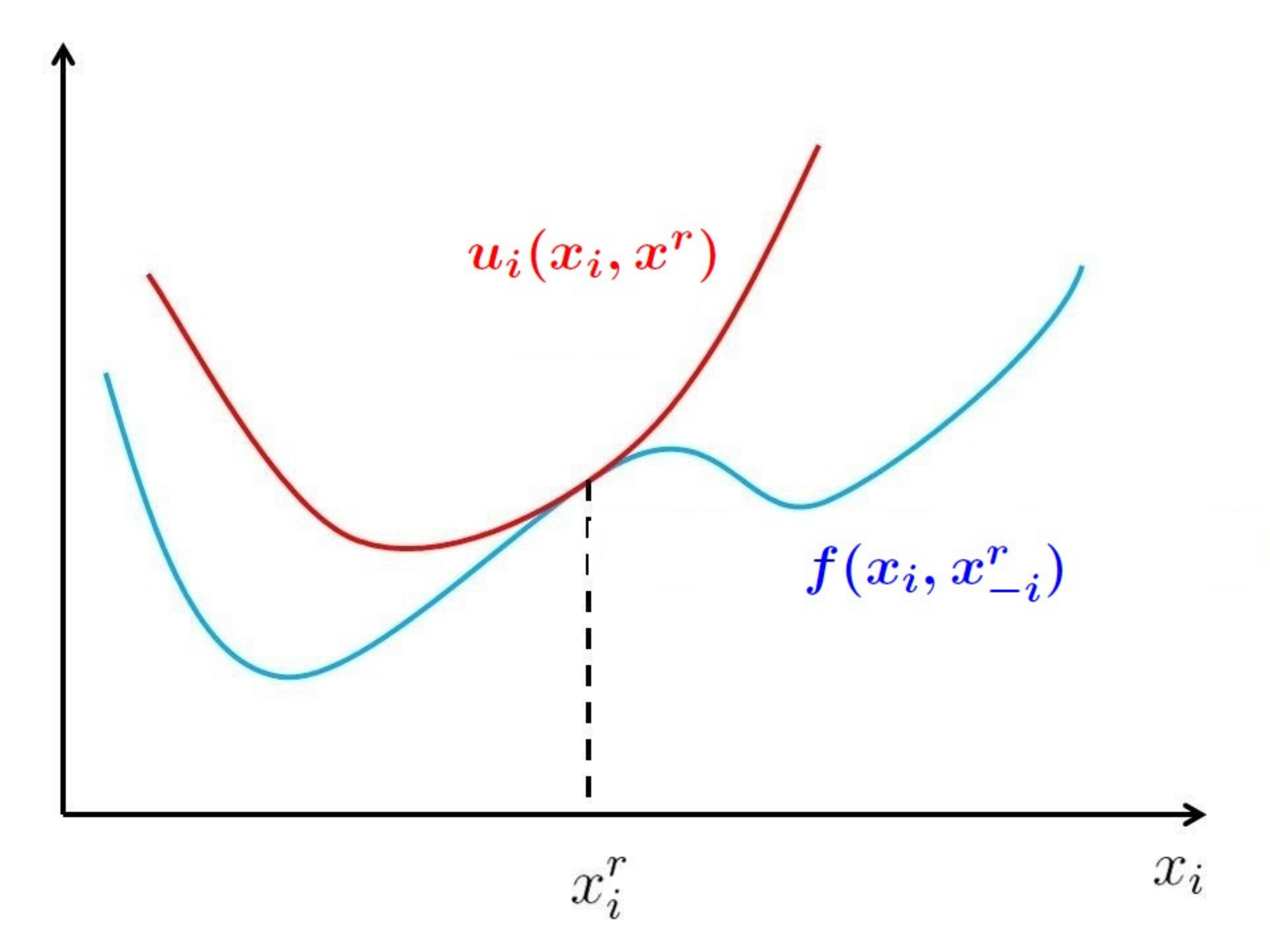}}}}
\hspace{1pc} { \subfigure[]
[]{\resizebox{.3\textwidth}{!}{\includegraphics{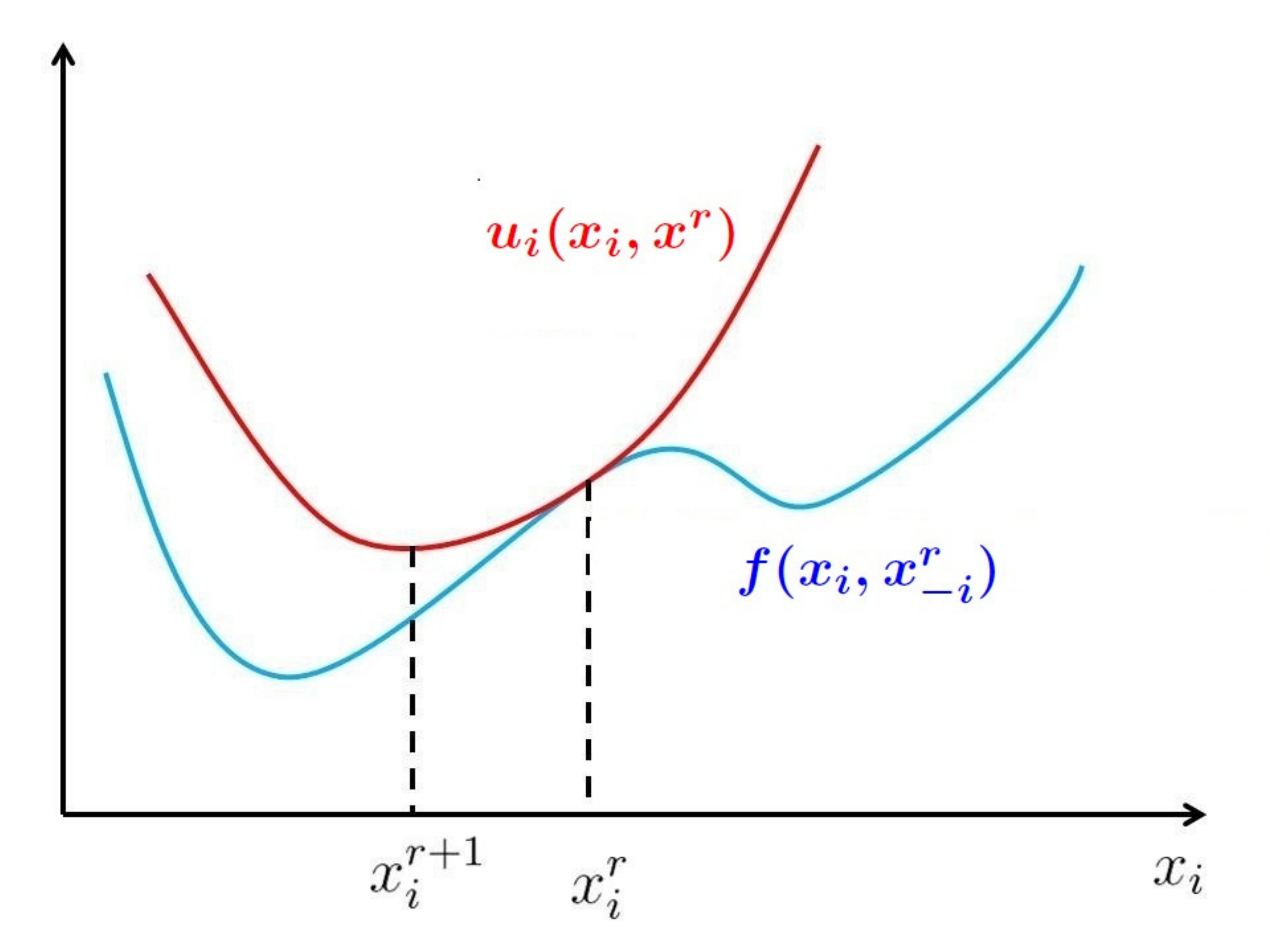}}}}\hspace{1pc}
\end{center}
\caption{\footnotesize Illustration of the upper-bound minimization step of the BSUM method. Here we assume that coordinate $i$ is updated at iteration $r+1$. It is clear from the figure that after solving the BSUM subproblem \eqref{eq:BSUMperblock}, $f(x^{r+1}_i, x^r_{-i})<f(x^{r}_i, x^r_{-i})$, that is, the objective function is strictly decreased.}\label{fig:upperbound-step}
\end{figure}

Now that we have seen the main steps of the algorithm, next we describe its theoretical properties.  We address the following questions related to the convergence of the BSUM: When does the BSUM converge? How fast does it converge? When does the BSUM fail and why? The answers to these questions will be instrumental in understanding the framework as well as evaluating the performance of various related algorithms.}

\subsection{When does the BSUM converge?}\label{sub:BSUM:converge}

In order to discuss the convergence property of the algorithm, we first investigate the optimality conditions for problem \eqref{eq:OriginalProblem}, which characterize the set of solutions that we would like our algorithm to reach. To this end, we introduce two related notions, one is called the {\it stationary solution} and the other is the {\it coordinatewise minimum solution}; see Table \ref{table:def} for their precise definitions. Roughly speaking, the coordinatewise minimum $\hat{x}$ is the point where no single block $\hat{x}_i$, $i=1,\cdots, n$ has a better direction to move to, while at a stationary point $x^*$ the entire vector cannot move to a better direction. Further, a stationary solution must be a coordinatewise minimum, but the reverse direction is generally not true; see the example next.

\begin{exa}\label{ex:regularity}
 {\it Consider the convex function $f(z) = \|Az\|_1$, where $A = [3 \;4; 2\; 1] \in \mathbb{R}^{2 \times 2}$. Obviously this function has a unique stationary solution $(z_1, z_2)=(0,0)$, which is also the global optimal solution. Further, the point $\hat{z} = (-4,3)$ is a coordinatewise minimum with respect to the two standard coordinates since $f'(\hat{z};d)\geq 0, \;\forall d\in \{(d_1,d_2)\in \mathbb{R}^2 \mid d_1 d_2 = 0\}$; but this point is not a stationary solution as $f'(\hat{z};\hat{d})<0$ for $\hat{d} = (4,-3)$. This fact can be also observed in the contour plot of the function in Fig.~\ref{fig:Contour}.}
\end{exa}
 \begin{figure}
 \centering
 \vspace{-0.4cm}
\includegraphics[width=0.4\linewidth]{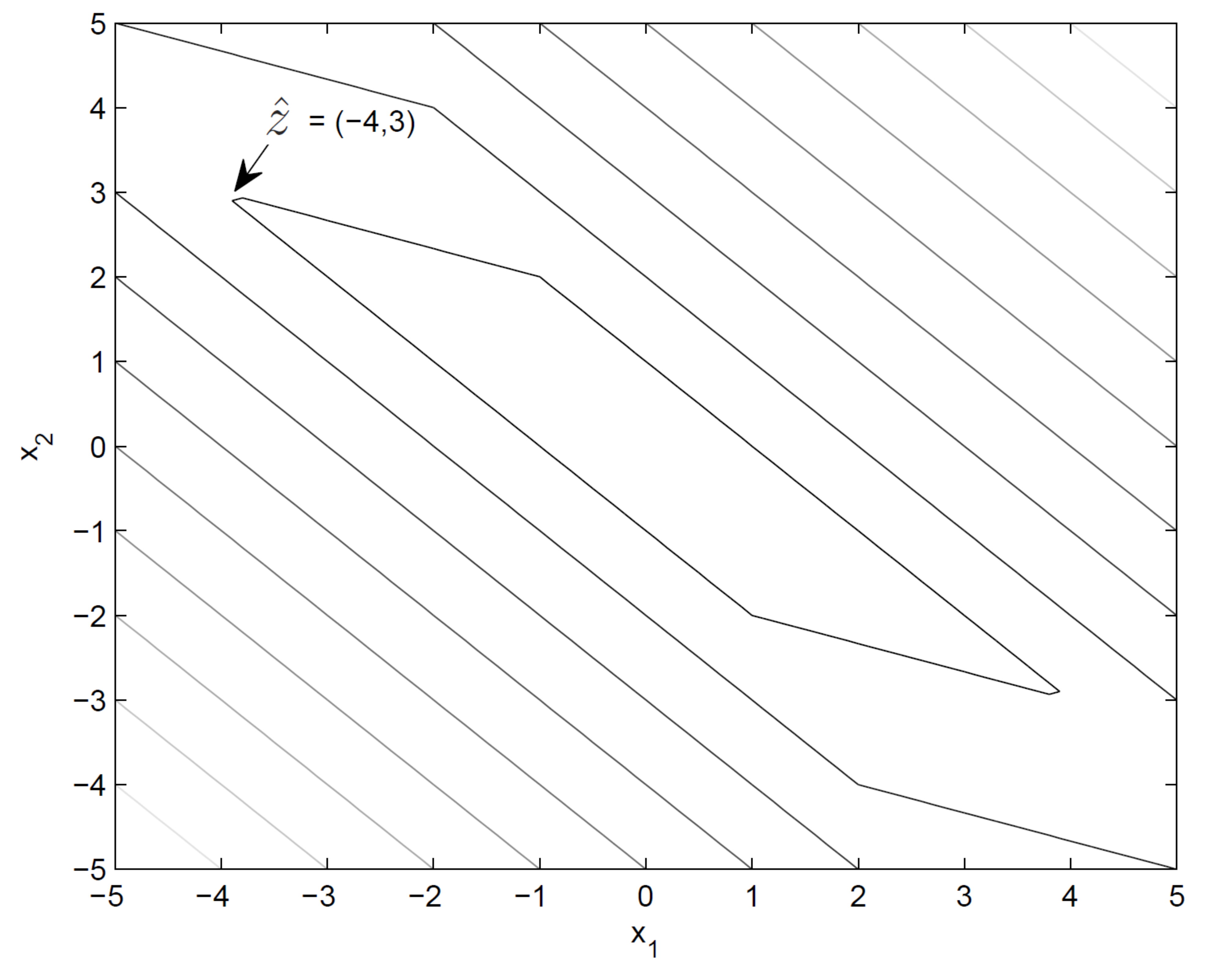}
\vspace*{-15pt}
\caption{ \footnotesize \cite{Razaviyayn12SUM} The contour plot of the function $f(z) = \|Az\|_1$,
where $A = [3 \;4; 2\; 1] \in \mathbb{R}^{2 \times 2}$. The point $\hat{z} =
(-4,3)$ is a coordinatewise minimum but not a stationary solution.\protect\footnotemark}\label{fig:Contour}
\end{figure}
\footnotetext{Copyright $\copyright$2013 Society for Industrial and Applied Mathematics.  Reprinted with permission.  All rights reserved.}

{This example confirms that the coordinatewise minimum  can be much inferior to the stationary solution.} Therefore in the subsequent discussion we will mainly focus on finding the stationary solutions rather than the coordinatewise minimum. An immediate question is that: whether one can easily distinguish these two types of solutions, or for that matter, when does a coordinatewise minimum become a stationary solution? {Let us define a {\it regular point} $x\in \mathcal{X}$ as a point that satisfies the following statement: if $x$ is coordinatewise minimum, then it is a stationary solution.} It turns out that for a large and popular subclass of problem \eqref{eq:OriginalProblem} expressed below in \eqref{eq:problemBSUM2} where the nonsmooth function is {\it separable} across the blocks, all feasible points $x\in\cX$ are regular.
\begin{align}
\min\; f(x):= g(x_1,\cdots, x_n)+\sum_{i=1}^{n}h_i(x_i),\quad\st\ x_i\in \cX_i, \ i=1,\cdots, n\label{eq:problemBSUM2}.
\end{align}
The main convergence result for BSUM method is given below, which is adapted from \cite[Theorem 2]{Razaviyayn12SUM}

\begin{thm}\label{thm:main}
 {\it Suppose the cyclic coordinate selection rule is chosen, i.e., $\cI^r = \{(r \;{\rm mod} \;n)+1\}$. Let $\{x^r\}_{r=1}^{\infty}$ be a sequence generated by the BSUM algorithm. Suppose Assumption A holds, and that each $x^r$ is regular. Then the following is true:}
\begin{itemize}
\item[(a)] {\it Suppose that the function $u_i(x_i,y)$ is quasi-convex in~$x_i$ for $i=1,\ldots,n$. Furthermore, assume that the subproblem~\eqref{eq:BSUMperblock} has a unique solution for any point~$x^{r-1} \in \mathcal{X}$. Then every limit point~$x^*$ of $\{x^r\}$ is a stationary point of \eqref{eq:OriginalProblem}.}
\item[(b)] {\it Suppose the level set $\mathcal{X}^0 = \{x \mid f(x) \leq f(x^0)\}$ is compact. Furthermore, assume that the subproblem~\eqref{eq:BSUMperblock} has a unique solution for any point $x^{r-1} \in \mathcal{X}$, $r\ge 1$ for at least $n-1$ blocks. Then $\{x^r\}$ converge to the set of stationary points, i.e.,
\[
{
\lim_{r \rightarrow \infty} \quad d(x^r,\cX^*) = 0,}
\]}
{\it where $d(x,\cX^*) \triangleq \min_{x^* \in \cX^*} \|x - x^*\|$  and $\cX^*$ is the set of stationary points. }
\end{itemize}
\end{thm}
Here the first part of the result deals with the possibility of an unbounded sequence, whereas the second part assumes that the sequence lies in a compact set, therefore it has a slightly stronger claim. {Note that Theorem \ref{thm:main} can be easily extended to all the coordinate selection rules given in Table \ref{table:coordinate}, and for the randomized version the convergence is with probability 1. }

The readers should pay special attention to the conditions set forth by Theorem~\ref{thm:main}. First it says that the upper bound function needs to be picked carefully to satisfy both Assumption A and the requirement that at least $n-1$ subproblems \eqref{eq:BSUMperblock} have a unique solution. However, when the objective function $f(x)$ is convex, the uniqueness requirement of the per-block solution can be dropped; see \cite{tseng01}. Also Theorem~\ref{thm:main} requires the problem to be well-defined so that coordinatewise optimal solutions are equivalent to the stationary solution (which precludes objective function $\|A x\|_1$  in Example \ref{ex:regularity}). In the next subsection we will demonstrate, through a couple of concrete examples, that relaxing some of these conditions indeed leads to the divergence of the BSUM.

\subsection{When does BSUM fail?}\label{sub:BSUM:fail}
In this subsection, we provide a few examples in which BSUM fails to converge to any stationary solutions. The examples in this section serve as reminders that in practice, in order to avoid those pitfalls, extreme care must be exercised when designing large-scale optimization algorithms.

Our first example comes from Example \ref{ex:regularity}. It demonstrates the consequence of lacking the regularity condition.
\begin{exa}\label{ex:diverge:regularity}
{\it Consider the following unconstrained convex optimization problem $\min_{x}\; \left\|\sum_{i=1}^{2}A_i x_i\right\|_1$, where $A_1=[3\; 2]^T$ and $A_2 = [4\; 1]^T$. Clearly by defining $A=[3\; 4; 2\; 1]$, the objective function can be rewritten as $\|Ax\|_1$, which is not regular at point $(x_1, x_2)=(-4,3)$ (cf. Example \ref{ex:regularity}). It follows that by setting $u_i(x_i,z)=\left\|\sum_{i=1}^{2}A_i x_i\right\|_1$ (no approximation is used), and letting $(x^0_1, x^0_2)=(-4,3)$, the BSUM algorithm will not be able to move forward for either $x_1$ or $x_2$, thus it will be stuck at the non-interesting point $(-4,3)$, far away from the only stationary solution $(0,0)$.
}
\end{exa}

The next example shows that BSUM fails to converge if the feasible set $\cX$ is no longer a Cartesian product of feasible sets $\cX_1,\cdots, \cX_n$, a fact that we have taken for granted so far.

\begin{exa} \label{ex:diverge:coupling}{\it
Consider the following simple quadratic problem:
  \begin{align}
  \min \quad & x^2_1+x^2_2, \quad  {\st}\quad x_1+x_2=2.\nonumber
  \end{align}
The optimal objective value is $2$. Consider the BSUM algorithm with an arbitrary approximation function satisfying Assumption A, but initiated at the point $(x^0_1,x^0_2)=(0,2)$. Then the BSUM method will be stuck at this non-interesting point without making any progress, because it is not possible to change a single block without violating the coupling constraint. }
\end{exa}

{Our next example shows that BSUM could diverge if the approximation function $u_i$ violates Assumption A.
\begin{exa} \label{ex:diverge:upper bound}{\it
Consider the simple quadratic problem   \begin{align}
  \min \quad & x^2_1+x^2_2 + 2x_1 x_2, \quad  {\st}\quad -1\le x_1, x_2\le 1. \nonumber
  \end{align}
 The optimal objective value is $0$, with $(x^*_1, x^*_2)=(0,0)$. Consider using the BSUM algorithm with a {\it linear} approximation function, which violates Assumption (A2). More specifically, for a given feasible tuple $(x_1,x_2)=(\hx_1,\hx_2)$, the $x_1$'s subproblem becomes
\begin{align*}
\min_{x_1} \quad & \langle \hx_1+ \hx_2, x_1\rangle, \quad {\st} \quad -1\le x_1\le 1.
\end{align*}
Clearly the optimal solution is either $x_1=-1$ or $x_1=1$. The same happens when solving the subproblem for $x_2$. Therefore the algorithm will never reach the optimal solution $(x^*_1, x^*_2)=(0,0)$. Further, if the feasible sets of $x_1$ and $x_2$ are unbounded, then the BSUM can diverge.
}
\end{exa}}

Our last example is taken from Powell \cite{PowellBCD73}. It shows that without the uniqueness assumption of each subproblem \eqref{eq:BSUMperblock}, the BSUM algorithm is not necessarily convergent.
\begin{exa} \label{ex:diverge:unique}{\it
Consider the following unconstrained problem
\begin{align}
\min\quad f(x)&:= -x_1 x_2-x_2x_3 - x_3x_1 + (x_1-1)^2_{+}+(-x_1-1)^2_{+}\nonumber\\
&\quad +(x_2-1)^2_{+}+(-x_2-1)^2_{+}+(x_3-1)^2_{+}+(-x_3-1)^2_{+}\nonumber
\end{align}
where the notation $(z)^2_{+}$ means $(\max\{0, z\})^2$. In this case, fixing $(x_2, x_3)$ and optimizing over $x_1$ yields the following solution
\begin{align}\label{eq:solution:x}
x_1 = \left\{ \begin{array}{ll}(1+\frac{1}{2}|x_2+x_3|)\mbox{\rm sign}(x_2+x_3), & \mbox{if}\; x_2+x_3\ne 0,\\
\rm{[-1,  \;1]}, &\mbox{otherwise}
\end{array}\right..
\end{align}
Similar solution can be obtained for $x_2$ and $x_3$ as well. Suppose we set $u_i(x_i,z)=f(x)$ for all $i$ (no approximation is used), and letting $(x^0_1, x^0_2, x^0_3)=(-1-\epsilon, 1+\frac{1}{2}\epsilon, -1-\frac{1}{4}\epsilon)$ for some $\epsilon>0$. Then it can be shown the applying the {cyclic version of the} BSUM algorithm, the iterates will be cycling around six points $(1, 1, -1)$, $(1, -1, -1)$, $(1, -1, 1)$, $(-1, -1, 1)$, $(-1, 1, 1)$, $(-1, 1, -1)$, and none of these six points is a stationary solution of the original problem. The reason for such pathological behavior is that, in any one of the six limit points above, there are at least two subproblems that have multiple optimal solutions. For example, at $(1, 1, -1)$ and fixing $x_2,x_3$ (resp. $x_1,x_3$), the optimal solution for $x_1$  (resp. $x_2$) is any element in the interval $[-1,\;1]$; cf. \eqref{eq:solution:x}.}
\end{exa}

A natural question at this point is, can we make the BSUM work for these examples? The answer is affirmative, but how this can be done requires a case by case study. To handle the first two examples (i.e., Example \ref{ex:diverge:regularity} -- \ref{ex:diverge:coupling}), a generalized version of BSUM is needed, which will be discussed in section \ref{sec:extension}. For the third example, one can simply pick a better upper bound to guarantee convergence. For example, if we pick the proximal upper bound (cf. Table \ref{table:bounds}):  $u_i(x_i,z)=f(x)+\frac{\gamma}{2}\|x_i-z_i\|^2$, then each subproblem will have a unique solution, and the algorithm will not be trapped by the non-interesting solutions given in Example \ref{ex:diverge:unique}. Notice that the use of $\frac{\gamma}{2}\|x_i-z_i\|^2$ inhibits the move of $x_i$ from its current position $z_i$. Thus, the main message here is that when optimizing each block, it is beneficial, at least theoretically, to be {\it less greedy} so that a conservative step is taken towards reducing the objective. The extent of the ``conservativeness" for the per block update is then naturally characterized by the chosen upper bounds.
{Quite interestingly, the difficulty with the non-unique subproblem solution can also be resolved by using randomization. Formally, we have the following corollary to Theorem \ref{thm:main}}.
\begin{coro}
{\cite{razaviyayn14thesis} Suppose the level set $\mathcal{X}^0 = \{x \mid f(x) \leq f(x^0)\}$ is compact. Then under the randomized block selection rule, the iterates generated by the BSUM algorithm converge to the set of stationary points almost surely, i.e.,
\[
\lim_{r\rightarrow \infty} d(x^r,\cX^*) = 0,\quad {\rm almost \;surely}.
\]}
\end{coro}

\subsection{How fast does the BSUM converge?}\label{sub:BSUM:converge:rate}

Now that we have examined the convergence of the BSUM, let us proceed next to characterize the convergence speed of the algorithm.  No doubt that this is an important issue, especially so for big data optimization -- the sheer size of the data and limited computational resource makes it difficult to optimize a problem to global optimality. Consequently, we are generally satisfied with high-quality suboptimal solutions, and are mostly concerned with how quickly such solutions can be obtained.

Recently, extensive research efforts have been devoted to analyzing the convergence rate for various BSUM-type algorithms, most of which use randomized coordinate selection rules and/or quadratic upper bound functions (cf. Table \ref{table:bounds}) to solve convex problems; for example see \cite{nestrov12,richtarik12,lu13complexity, Beck13,Beck13b,hong13complexity,Saha10}. Since it is not possible to go over all the details of these different variations of BSUM, here we present in a high level a fairly general result that covers a large family of upper bound functions satisfying Assumption A, as well as all coordinate selection rules outlined in Table \ref{table:coordinate}.

To this end, let us make the following additional assumptions.

\noindent{\bf Assumption B}
\begin{itemize}
\item[(B1)] $f(x):=g(x)+\sum_{i=1}^{n}h_i(x_i)$, where $g(x)$ is a smooth convex function with Lipschitz continuous gradient, i.e. there exists a constant $L$ such that $\|\nabla g(x)-\nabla g(y)\|\le L\|x-y\|, \; \forall~x,y\in \cX.$
    Further both $g$ and $h_i$'s are convex functions.
\item[(B2)] The level set $\{x\mid f(x)\le f(x^0), x\in \cX\}$ is compact.
\item[(B3)] Each upper bound function $u_i(x_i, z)$ is {\it strongly convex} with respect to $x_i$.
\end{itemize}

An {\it $\epsilon$-optimal solution} $x^{\epsilon}\in\cX$ is defined as an $x^{\epsilon}\in\{x\mid x\in \cX, f(x)-f(x^*)\le \epsilon\}$, where $f(x^*)$ is the globally optimal objective value of problem \eqref{eq:problemBSUM2}. Suppose both Assumptions A and B are satisfied. Then it is shown in \cite{hong13complexity, razaviyayn14thesis} that BSUM takes at most ${c}/{\epsilon}$ iterations to find an $\epsilon$-optimal solution, for all coordinate rules specified in Table \ref{table:coordinate}, where $c>0$ is a constant only related to the description of the problem. Such type of convergence rate is known as {\it sublinear convergence}. Here the main message is that under Assumptions A and B, the algorithm generally converges sublinearly in the order of $1/\epsilon$. Further, for different special forms of BSUM, the constant $c$ in front of $1/\epsilon$ can be significantly refined so that it is {\it independent} of problem dimension; see \cite{nestrov12,richtarik12}. It is also interesting to note here that when the objective $f$ is either strongly convex, or convex but with certain special structure, the BSUM  achieves {\it linear rate of convergence}. That is, BSUM takes at most $\mathcal{O}(\log(c/\epsilon))$ iterations to find an $\epsilon$-optimal solution, which is much faster than the sublinear rate; see e.g., \cite{hong13BSUMM,Sanjabi13} for the related discussions. 

Finally, we briefly mention that it is possible to relax certain conditions in Assumption B to obtain refined rates.  For example, \cite{hong13complexity} shows that dropping the per-block strong convexity assumption in (B3) still achieves an $\mathcal{O}({1}/{\epsilon})$ sublinear convergence. In \cite{razaviyayn14thesis,xu12} it is shown that when removing the convexity Assumption (B1), it is also possible to characterize the convergence rate to stationary solutions. {In \cite{hong13complexity} it is shown that when there are two blocks of variables, the cyclic version of the BSUM can be accelerated to achieve an improved $\mathcal{O}(1/\sqrt{\epsilon})$ complexity.} {In a few recent works \cite{lin14, Fercoq13}, it is shown that when randomized block selection and the quadratic upper bound are used, it is possible to accelerate the BSUM with any $n>2$ blocks. }


\section{Algorithms Covered by the BSUM Framework}\label{sec:algorithms}
Now that we have a good understanding of the theoretical properties of the family of BSUM algorithms, we demonstrate in this section the wide applicability of the BSUM framework by revealing its connection to a few well-known algorithms for high dimensional massive data analysis. For each of the examples presented below, we pay special attention to the design of the upper bound functions.

 \subsection{The Block Coordinate Descent (BCD) algorithm}
 The first algorithm that the BSUM covers is obviously the classic BCD described in Section \ref{sub:bcd}. Here the upper bound function is simply the original objective itself, i.e., $u_i(x_i,z):= f(x_i,z_{-i}), \; \forall~x_i\in\cX_i, z\in \cX, \forall~i$. We would like to mention that all the convergence and rate of convergence analysis of BSUM carries over to the BCD method. Specifically, the result in Section \ref{sub:BSUM:converge:rate} implies that the BCD method (with coordinate update rules specified in Table
 \ref{table:coordinate}) converges in a sublinear manner whenever Assumptions (B1) and (B2) are satisfied. This result by itself is interesting, as there has been limited theoretical analysis for general form of BCD algorithm when applied to problems satisfying Assumptions (B1) and (B2).

 \subsection{The Convex-Concave Procedure (CCCP)}
 Consider the following unconstrained nonconvex problem, known as the difference of convex (DC) program:
 $\min \; f(x):= g_1(x)-g_2(x)$ where both $g_1$ and $g_2$ are convex functions.  The well-known CCCP algorithm \cite{CCCP} generates a sequence $\{x^r\}$ by solving
\begin{equation*}
x^{r+1} = \arg \min_{x} \; u(x,x^r),
\end{equation*}
where $u(x,x^r) := g_{1}(x) - \langle x - x^r, \nabla g_{2}(x^r)\rangle - g_{2}(x^r)$. Clearly, the {\it linear upper bound} in Table \ref{table:bounds} is used here, therefore  CCCP is a special case of the BSUM algorithm with a single block of variables. Furthermore, the updates can also be done in a block coordinate manner. 

 \subsection{The majorization-minimization (MM) algorithm}
 The MM algorithm which has been widely used in statistics \cite{Hunter04}, can be viewed as the single block version of BSUM. Consider the problem of $\min_{x\in\cX} f(x)$ where $f(x)$ is a smooth function. The basic idea of the MM algorithm is to first construct a ``majorization" function $u(x,z)$ for the original objective $f(z)$, such that
 \begin{align}\label{eq:MM}
 u(x,z)\ge f(z), \forall~x, z\in \cX, \quad u(x,x)=f(x), \; \forall~x\in \cX.
 \end{align}
 Then the algorithm generates a sequence of iterates by successively minimizing $u(x,x^r)$. Clearly this algorithm represents a slight generalization of the CCCP discussed in the previous subsection, but nevertheless still falls in the framework of BSUM.

 As another concrete example of the MM algorithm, let us consider the celebrated expectation maximization (EM) algorithm \cite{EMDempster}.  Let $w$ be a vector observation used for estimating the parameter $x$. The maximum likelihood estimate of $x$ is given as (where $p(w|x)$ denotes the conditional probability of $w$ given $x$)
\begin{equation}\label{EQ:EMML}
\hat{x}_{\rm ML} = \arg \max_{x} \; \ln p(w|x).
\end{equation}
Let the random vector $z$ be some hidden/unobserved variable. The EM algorithm generates a sequence $\{x^r\}$ by iteratively performing the following two steps
{\it (1)}{ E-Step: Calculate $g(x,x^r) := \mathbb{E}_{z|w,x^r} \{\ln p(w,z|x)\}$; {\it(2)} M-Step: $x^{r+1} = \arg \max_x \; g(x,x^r)$.
To see why the EM algorithm is a special case of MM (hence a special case of BSUM) let us rewrite problem \eqref{EQ:EMML} as:
$\min_x \; -\ln p(w|x)$, whose objective can be bounded by
\begin{align}
-\ln p(w|x) & = -\ln \;\;\mathbb{E}_{z|x}\; p(w|z,x)\nonumber\\
& = -\ln \;\;\mathbb{E}_{z|x} \left[\frac{p(z|w,x^r) p(w|z,x)}{p(z|w,x^r)}\right]\nonumber\\
& = -\ln \;\;\mathbb{E}_{z|w,x^r} \left[\frac{p(z|x) p(w|z,x)}{p(z|w,x^r)}\right]\nonumber\\
& \leq - \mathbb{E}_{z|w,x^r} \ln\left[\frac{p(z|x)p(w|z,x)}{p(z|w,x^r)}\right]\nonumber\\
& = - \mathbb{E}_{z|w,x^r} \ln p(w,z|x) + \mathbb{E}_{z|w,x^r} \ln p(z|w,x^r)\nonumber\\
& := u(x,x^r),\nonumber
\end{align}
{\color{black}where the inequality is due to the fact that a convex function $f$ must satisfy $\mathbb{E}[f(x)]\ge f(\mathbb{E}[x])$ (by the Jensen's inequality)}. Since $\mathbb{E}_{z|w,x^r} \ln p(z|w,x^r)$ is not a function of $x$, the M-step can be written as
\[
x^{r+1} = \arg\min_{x} u(x,x^r).
\]
Furthermore, it is not hard to see that $u(x^r,x^r) = - \ln p(w|x^r)$, therefore, both conditions in \eqref{eq:MM} are satisfied. Similarly as in the previous case, one can extend the EM to a block coordinate version; see \cite{Razaviyayn12SUM} for detailed discussion.

 \subsection{The Proximal Point Algorithm (PPA)}
 The classical PPA (see, e.g.,
\cite[Section 3.4.3]{Bertsekas_Book_Distr}) obtains a solution of
the problem $\min_{x\in\mathcal{X}}f(x)$ by solving the following equivalent problem
\begin{align}
\min_{x\in\mathcal{X},y\in\mathcal{X}}f(x)+\frac{\gamma}{2}\|x-y\|_2^2,\label{problemProximal}
\end{align}
where $f(\cdot)$ is a convex function, $\mathcal{X}$ is a closed
convex set, and $\gamma>0$ is a coefficient. Clearly problem
\eqref{problemProximal} is strongly convex
in both $x$ and $y$ so long as $f(x)$ is convex (but not jointly strongly convex in $(x,y)$).
This problem can be solved by performing the following two steps alternatingly
\begin{align}
x^{r}&=\arg\min_{x\in\mathcal{X}}\left\{f(x)+\frac{1}{2c}\|x-y^{r-1}\|^2_2\right\}, \quad y^{r}=x^{r+1}.
\end{align}
Equivalently, the PPA algorithm can be viewed as successively minimizing the single block version of the proximal upper bound $u(x; x^r)$ given in Table \ref{table:bounds}. {Note that for a problem with a single block of variables, if $x\in \mathcal{X}$ is coordinatewise minimum, then it must be a global minimum solution.  Therefore every feasible solution $x\in\mathcal{X}$ is regular, and the convergence of PPA is covered by BSUM.}

Further, the PPA can be generalized to solve the multi-block problem \eqref{eq:OriginalProblem},
where $f(\cdot)$ is convex in each of its block components, but not
necessarily strictly convex. Directly applying BCD may fail to find a stationary solution for this
problem, as the per-block subproblems can contain multiple
solutions (cf. Example \ref{ex:diverge:unique}). Alternatively, we can consider an {\it alternating
PPA} \cite{ProximalBCD}, which successively solves the following subproblem
\begin{align*}
\min_{x_i}&\quad f(x_i,x_{-i}^r)+\frac{\gamma}{2}\|x_i-x^r_i\|^2_2,\quad
{\rm s.t.}\quad x_i\in\mathcal{X}_i\nonumber.
\end{align*}
Clearly this algorithm is a special form of BSUM with strongly convex proximal upper bound (cf. Table \ref{table:bounds}). It follows that each subproblem has a unique optimal solution, and by Theorem \ref{thm:main} it must converge to a stationary solution.

 \subsection{The Forward-Backward Splitting (FBS) Algorithm}
 The FBS algorithm (also known as the proximal splitting algorithm, see, e.g., \cite{Combettes09} and the references therein) for nonsmooth optimization solves the composite problem \eqref{eq:problemBSUM1} with a single block of variables (i.e., $n=1$),
where $h$
is convex and lower semicontinuous; $g$ is convex and has
Lipschitz continuous gradient, i.e., $\|\nabla g(x)-\nabla
g(y)\|\le L \|x-y\|$, $\forall~x,y\in\mathcal{X}$
and for some $L>0$.

Define the proximity operator ${\rm
prox}_{h}:\mathcal{X}\to\mathcal{X}$ as
\begin{align}\nonumber
{\rm
prox}_{h}(\bx)=\arg\min_{y\in\mathcal{X}} \; h(y)+\frac{1}{2}\|x-y\|^2.
\end{align}
The FBS iteration is given below
\cite{Combettes09}:
\begin{align}
x^{r+1}=\underbrace{{\rm prox}_{\beta h}}_{\rm backward~
step}\underbrace{(x^r-\beta\nabla g(x^r))}_{\rm
forward~step}\label{eq:forwardbackward}
\end{align}
where $\beta\in(0,1/L]$. Define 
\begin{align}
u(x,x^r):=
h(x)+\frac{1}{2\beta}\|x-x^r\|^2+\langle x-x^r,\nabla
g(x^r)\rangle+g(x^r),
\end{align}
which is the quadratic upper bound in Table \ref{table:bounds}, with $\bfPhi_1:=\frac{1}{\beta}\mathbf{I}$.
It is easy to see that the iteration \eqref{eq:forwardbackward} is
equivalent to the following iteration $
x^{r+1}=\arg\min_{x\in\mathcal{X}}u(x,x^r),
$ therefore it again falls under the BSUM framework.

Similar to the previous example, we can generalize the
FBS algorithm to solve multiple block problems of the form \eqref{eq:problemBSUM2}. 
The resulting algorithm, sometimes also known as the block coordinate proximal gradient (BCPG) method, has received significant attention recently due to its efficiency for solving certain big data optimization problem such as LASSO \cite{Friedman10}. We refer the readers to \cite{nestrov12,lu13complexity,lu13randomized} for recent developments and applications for BCPG.

Here we make a special note that by appealing to the general convergence rate result in Section \ref{sub:BSUM:converge:rate}, the BCPG method with any coordinate selection rules in Table \ref{table:coordinate} gives a sublinear convergence rate, when it is used  to solve problem \eqref{eq:problemBSUM2} that satisfies Assumption B.

 \subsection{The Nonnegative Matrix Factorization (NMF) Algorithm}\label{sub:nmf}
 Consider the following NMF problem:
 \begin{align}\label{eq:nmf}
 \min_{W\in\mathbb{R}^{M\times K}, H\in\mathbb{R}^{K\times M}}f(W, H):=\frac{1}{2}\|V-WH\|_F^2, \quad \st\; W\ge 0, \; H\ge 0,
 \end{align}
 where $V\in\mathbb{R}^{M\times N}_{+}$ is given. The problem has been extensively studied since Lee and Seung's seminal work \cite{Lee00algorithmsfor}, and it has wide applications in factor analysis, dictionary learning, speech analysis and so on \cite{mairal2009online}. In \cite{Lee00algorithmsfor}, a simple and efficient multiplicative algorithm is proposed:
 \begin{subequations}
 \begin{align}
 [H^{r+1}]_{j,i} &= [H^r]_{j,i}\frac{[(W^r)^T V]_{j,i}}{[(W^r)^T W^r H^r]_{j,i}}, \quad j=1,\cdots, K,\; i=1, \cdots, N\label{eq:Hupdate}\\
 [W^{r+1}]_{j,i} &= [W^r]_{j,i}\frac{[V(H^{r+1})^T ]_{j,i}}{[W^r H^{r+1} (H^{r+1})^T]_{j,i}}, \quad j=1,\cdots, M,\; i=1, \cdots, K\label{eq:Wupdate}.
 \end{align}
 \end{subequations}
 Here $[H^{r+1}]_{j,i}$ means the $(j,i)$th component of matrix $H^{r+1}$.
Below we show that when the iterates are well-defined (i.e., $[W^r_{ij}]>0$ and $[H^r_{ij}]>0$), the NMF iteration \eqref{eq:Hupdate} -- \eqref{eq:Wupdate} is also covered by BSUM \cite{Li12NMFnote}.

Let $H_i$ and $V_i$ represent the $i$th column of $H$ and $V$, respectively. Then at a given iterate $\{W^r, H^r\}$, the subproblem for optimizing $H_i$ is given by
\begin{align}
\min_{H_i\ge 0}\;f(H_i,\{W^r, H^r\}):=\frac{1}{2}\|V_i - W^r H_i\|_F^2+\frac{1}{2}\sum_{j\ne i}\|V_j - W^r H^r_j\|_F^2\label{eq:nmf:block}.
\end{align}

Define the upper bound function $u_i(H_i,\{W^r, H^r\})$ as
\begin{align*}
u_i(H_i,\{W^r, H^r\})&:=
f(H^r_i,\{W^r, H^r\})+(H_i-H^r_i)^T\nabla_{H_i} f(H_i,\{W^r, H^r\})\nonumber\\
&\quad+\frac{1}{2}(H_i-H^r_i)^T \bfPhi_i(W^r, H^r_i)(H_i-H^r_i)
\end{align*}
where $\bfPhi_i(W^r, H^r_i)$ is a diagonal matrix given by
\begin{align*}
\bfPhi_i(W^r, H^r_i): =\mbox{Diag}\left(\frac{[(W^r)^T W^r H_i^r]_1}{[H^r_i]_1},\cdots,\frac{[(W^r)^T W^r H_i^r]_K}{[H^r_i]_K}\right).
\end{align*}

Clearly, $\bfPhi_i(W^r, H^r_i)\succ 0$, and it is easy to show that  $\bfPhi_i(W^r, H^r_i)\succ (W^r)^T W^r$, where $(W^r)^T W^r$ is the Hessian of the objective of \eqref{eq:nmf:block} \cite{Lee00algorithmsfor}.
This implies that $u_i(H_i,\{W^r, H^r\})$ is the {quadratic upper bound given in Table \ref{table:bounds} of $f(H_i,\{W^r, H^r\})$. Further, one can check that the subproblem that minimizes $u_i(H_i,\{W^r, H^r\})$ has a unique solution, given by \eqref{eq:Hupdate}.}
Similar analysis can be established for the $W$-block update rule as well. Therefore we conclude that the iterates \eqref{eq:Hupdate} -- \eqref{eq:Wupdate} are a special case of BSUM. Finally we note that it is also possible to use different upper bound functions to derive more efficient update rules for the NMF problem \eqref{eq:nmf}; for example see \cite{fevotte11} where both the concave upper bound and the Jensen's upper bound (cf. Table \ref{table:bounds}) are used.

 \subsection{The Iterative Reweighted Least Squares  (IRLS) Method}

 The IRLS method is a popular algorithm used for solving big data problems such as sparse recovery \cite{Daubechies10}. Consider the following problem
 \begin{align}
\min_{x}&\quad  h(x) + \sum_{j=1}^{\ell}\|A_j x + b_j\|_2,\quad \st \quad x\in \mathbb{R}^m \label{eq:irls}
 \end{align}
 where $A_j\in\mathbb{R}^{k_i\times m}$ and $b_j\in\mathbb{R}^{k_i}$, and $h(x)$ is some convex function not necessarily smooth. For a set of applications for this model, see \cite[Section 4]{Beck13b}. Consider the following {\it smooth approximation} of problem \eqref{eq:irls}:
\begin{align}
\min_{x}&\quad h(x)+g(x): = h(x) + \sum_{j=1}^{\ell}\sqrt{\|A_j x + b_j\|^2_2+\eta^2},\quad \st \quad  x\in \mathbb{R}^m \label{eq:irls:smooth}
 \end{align}
 where $\eta$ is some small constant and $g(x)$ denotes the smooth part of the objective.  The IRLS algorithm solves problem \eqref{eq:irls:smooth} by performing the following iteration
 \begin{align*}
 x^{r+1}=\arg\min_{x\in\mathbb{R}^{m}}\left\{h(x)+\frac{1}{2}\sum_{j=1}^{\ell}\frac{\|A_j x+b_j\|^2+\eta^2}{\sqrt{\|A_j x^r+b_j\|^2+\eta^2}}\right\}.
 \end{align*}
 Define the following function for $g(x)$:
 \begin{align}
 u(x, x^r) = \frac{1}{2}\left(\sum_{j=1}^{\ell}\frac{\|A_j x+b_j\|^2+\eta^2}{\sqrt{\|A_j x^r+b_j\|^2+\eta^2}}+\sqrt{\|A_j x^r+b_j\|^2+\eta^2}\right).\label{bound:irls}
 \end{align}
 It is clear that $g(x^r) = u(x^r, x^r)$, so Assumption (A1) is satisfied. To verify Assumption (A2), we apply the arithmetic-geometric inequality, and have
 \begin{align*}
 u(x, x^r) & = \frac{1}{2}\left(\sum_{j=1}^{\ell}\frac{\|A_j x+b_j\|^2+\eta^2}{\sqrt{\|A_j x^r+b_j\|^2+\eta^2}}+\sqrt{\|A_j x^r+b_j\|^2+\eta^2}\right)\nonumber\\
 &\ge \sum_{j=1}^{\ell}\sqrt{\|A_j x+b_j\|^2+\eta^2} = g(x), \; \forall~x\in\mathbb{R}^{m}.
 \end{align*}
 Then according to Example \ref{ex:separability}, Assumption (A3) is automatically true, therefore we have verified that $u(x,x^r)$ defined in \eqref{bound:irls} is indeed an upper bound function for the smooth function $g(x)$. It follows that the IRLS algorithm} corresponds to a single-block BSUM algorithm. Notice that using the BSUM framework we can easily generalize the IRLS to the multi-block scenario.


\section{Applications of the
BSUM Framework}\label{sec:application}
In this subsection, we briefly review a few applications of the BSUM framework in wireless communication, bioinformatics, signal processing, and machine learning.
\subsection{Wireless communication and transceiver design}
\label{sec:WMMSE}
Consider a multi-input multi-output interference channel with $K$ transmitter-receiver pairs. Let $M$  (resp. $N$) be the number of antennas at each transmitter (resp. receiver) and each transmitter $k$, $k=1,2,\ldots,K,$ is interested in transmitting one data stream to its own receiver. Let $\mathbf{x}_k \in \mathbb{C}^M$ be the transmitted signal of user~$k$; assuming linear channel model, the received signal of user~$k$ can be written as
\[
\mathbf{y}_k = \underbrace{\mathbf{H}_{kk} \mathbf{x}_k} _{\rm desired \; signal} + \underbrace{\sum_{j \neq k} \mathbf{H}_{kj} \mathbf{x}_j}_{\rm multi-user \; interference} + \underbrace{\mathbf{n}_k}_{ \rm noise},
\]
where $\mathbf{H}_{kj} \in \mathbb{C}^{N \times M}$ is the channel from transmitter $j$ to receiver $k$ and $\mathbf{n}_k \in \mathbb{C}^N$ denotes the additive white Gaussian noise at the receiver $k$ with distribution $\mathcal{CN} (\mathbf{0}, \sigma^2 \mathbf{I})$.

When linear beamformers are employed at the transmitters and receivers, the transmitted signal and the estimated received data stream can be respectively written as
\[
\mathbf{x}_k = \mathbf{v}_k s_k \quad {\rm and} \quad \widehat{s}_k = \mathbf{u}_k^H \mathbf{y}_k,
\]
where $\mathbf{v}_k \in \mathbb{C}^M$  and $\mathbf{u}_k \in \mathbb{C}^N$ are respectively the transmit and receive beamformers. Here the transmitted datastream and the estimated datastream  at the receiver are denoted by $s_k \in \mathbb{C}$ and $\widehat{s}_k \in \mathbb{C}$ respectively.

A crucial task in modern wireless networks is to design the transmit and receive beamformers $\mathbf{v}_k$ and $\mathbf{u}_k$ in order to maximize a given utility of the system. Here, for simplicity of presentation, we consider the sum rate utility function as our objective. Therefore, our goal is to solve the following optimization problem:
\begin{equation}
\label{eq:SumRateMaximization}
\begin{split}
\max_{\mathbf{u},\mathbf{v}} \;\; &\sum_{k=1}^K R_k (\mathbf{u},\mathbf{v}) \\
{\rm s.t.} \quad & \|\mathbf{v}_k\|^2 \leq P_k, \quad \forall k=1,2,\ldots,K,
\end{split}
\end{equation}
where $P_k$ is the total power budget of user $k$ and $R_k (\mathbf{u},\mathbf{v})$, which is the communication rate of user $k$, is given by
\[
R_k(\mathbf{u},\mathbf{v}) = \log\left( 1+ \frac{|\mathbf{u}_k^H \mathbf{H}_{kk} \mathbf{v}_k|^2}{\sigma^2 \|\mathbf{u}_k\|^2 + \sum_{j \neq k} |\mathbf{u}_k^H \mathbf{H}_{kj} \mathbf{v}_j|^2} \right).
\]
Problem \eqref{eq:SumRateMaximization} is nonconvex and known to be NP-hard \cite{luo08a}. Using the well-known relation between the signal to interference plus noise ratio (SINR) and the mean square error (MSE) value, one can rewrite \eqref{eq:SumRateMaximization} as \cite{shi11WMMSE_TSP,Baligh14}:
\begin{equation}
\begin{split}
\min_{\mathbf{v},\mathbf{u}} \;\; &\sum_{k=1}^K \log \left(e_k (\mathbf{u},\mathbf{v})\right) \\
{\rm s.t.} \quad & \|\mathbf{v}_k\|^2\leq P_k,\quad \forall k= 1,2,\ldots,K,
\end{split}
\end{equation}
where $e_k(\mathbf{u},\mathbf{v})$ is the MSE value and is given by
\[
e_k(\mathbf{u},\mathbf{v})  = |\mathbf{u}_k^H \mathbf{H}_{kk} \mathbf{v}_k -1 |^2 + \sum_{j\neq k} |\mathbf{u}_k^H \mathbf{H}_{kj} \mathbf{v}_j|^2 + \sigma^2.
\]
Since the $\log(\cdot)$ function is concave, it is upper bounded by its first order approximation (i.e., the linear upper bound in Table~\ref{table:bounds}). Therefore, we can define the function
\begin{equation}
\label{eq:WMMSEUpperBound}
u(\mathbf{x},\mathbf{x}^r) = \sum_{k=1}^K \left(\log\left( e_k (\mathbf{x}^r)\right)  + \left(e_k(\mathbf{x}^r)\right)^{-1} \left( e_k(\mathbf{x}) - e_k(\mathbf{x}^r)\right)\right)
\end{equation}
{where $\mathbf{x} \triangleq (\mathbf{u},\mathbf{v})$ is the optimization variable and $\mathbf{x}^r \triangleq (\mathbf{u}^r,\mathbf{v}^r)$ denotes the beamformer at iteration $r$. It is not hard to see that the approximation function in \eqref{eq:WMMSEUpperBound} is a valid upper-bound in  the BSUM framework  and at each iteration~$r$, this choice of approximation function leads to a quadratic programming problem which has closed form solutions. The resulting algorithm, dubbed weighted minimization of mean square error (WMMSE), converges to a stationary point of the problem and in practice  it typically converges in a few iterations \cite{shi11WMMSE_TSP} even for large size problems~\cite{Razaviyayn11grouping}.

The  reader is referred to \cite{shi11WMMSE_TSP,Baligh14,christensen08,hong12sparse,razaviyayn2013stochastic} for more details of the algorithm and its extensions to various beamformer design scenarios and different utility functions. It is also worth noting that many other interesting transceiver design algorithms also fall into the BSUM framework; see  \cite{shi2008pricingmiso,Shi:2009,kim11,Ng10} for more details.
}

\subsection{ Bioinformatics and Signal Processing}
Here we briefly outline two interesting big data applications of the BSUM framework in bioinformatics and signal processing.
\subsubsection{\color{black}Abundance Estimation in Modern High-Throughput Sequencing Technologies}
An essential step in the analysis of modern high throughput sequencing technologies of biological data is to estimate the abundance level of each transcript in the experiment. Mathematically,   this problem can be stated as follows. Consider $M$ transcript sequences $s_1,\ldots,s_M \in\{A,C,G,T\}^L$ with the corresponding abundance levels $\rho_1,\ldots,\rho_M$ such that $\sum_{m=1}^M \rho_m = 1$. Let $R_1,\ldots,R_N$ be noisy sequencing reads originated from the transcript sequences, where each read $R_n$, $n=1,\ldots,N$, is originated from only one of the transcript sequences $s_1,\ldots,s_M$. Given the observed reads, the likelihood of the abundance levels $\rho_1,\ldots,\rho_M$ can be written as
\begin{align}
{\rm Pr}\left(R_1,\ldots,R_N ; \rho_1,\ldots,\rho_M\right) & = \prod_{n=1}^N {\rm Pr}\left(R_n ; \rho_1\ldots \rho_M\right) \nonumber\\
&= \prod_{n=1}^N \left( \sum_{m=1}^M {\rm Pr}\left(R_n \mid {\rm read\;} R_n {\; \rm from \; sequence \;} s_m\right) {\rm Pr}(s_m)\right) \nonumber\\
&= \prod_{n=1}^N \left( \sum_{m=1}^M \alpha_{nm} \rho_m\right), \nonumber
\end{align}
where $\alpha_{nm} \triangleq {\rm Pr}\left(R_n \mid {\rm read\;} R_n {\; \rm from \; sequence \;} s_m\right)$  can be obtained efficiently using an alignment algorithm such as the ones based on the Burrows-Wheeler transform; see,  e.g., \cite{langmead2009ultrafast,li2009fast}. Therefore, given $\{\alpha_{nm}\}_{n,m}$, the maximum likelihood estimation of the abundance levels can be stated as
\begin{equation}
\label{eq:MaximumLikelihoodAbundance}
\begin{split}
\widehat{\rho}_{ML} =   \arg\min_{\rho}  \quad & -\sum_{n=1}^N \log\left( \sum_{m=1}^M    \alpha_{nm} \rho_m\right) \\
{\rm s.t.} \quad &\sum_{m=1}^M \rho_m = 1, \quad {\rm and }\quad \rho_m\geq 0, \;\forall m =1,\ldots,M.
\end{split}
\end{equation}
As a special case of the EM algorithm, a popular approach for solving this optimization problem is to successively minimize a  local tight upper-bound of the objective function. In particular, the eXpress software \cite{roberts2013streaming} solves the following optimization problem at the $r$-th iteration of the algorithm:
\begin{equation}
\label{eq:MaximumLikelihoodAbundanceUpperBd}
\begin{split}
{\rho}^{r+1}=   \arg\min_{\rho}  \quad & -\sum_{n=1}^N \left( \sum_{m=1}^M \left(   \frac{\alpha_{nm} \rho_m^r}{\sum_{m^\prime = 1}^M \alpha_{n m^\prime} \rho_{m^\prime}^r } \log \left(\frac{\rho_m}{\rho_m^r}\right) \right) + \log\left(\sum_{m=1}^M \alpha_{nm} \rho_m^r\right) \right)\\
{\rm s.t.} \quad &\sum_{m=1}^M \rho_m = 1, \quad {\rm and }\quad \rho_m\geq 0, \;\forall m =1,\ldots,M.
\end{split}
\end{equation}
Using Jensen's inequality, it is not hard to check that \eqref{eq:MaximumLikelihoodAbundanceUpperBd} is a valid upper-bound of \eqref{eq:MaximumLikelihoodAbundance} in the BSUM framework. Moreover, \eqref{eq:MaximumLikelihoodAbundanceUpperBd} has a closed form solution obtained by
\[
\rho_m^{r+1} = \frac{1}{N}\sum_{n=1}^N\frac{ \alpha_{nm} \rho_m^r}{ \sum_{m^\prime=1}^M  \alpha_{nm^\prime} \rho_{m^\prime}^r}, \quad  \forall m=1,\ldots, M,
\]
which makes the algorithm computationally efficient at each step.
{In practice, the above algorithm for abundance estimation converges in a few iterations.  Moreover, this algorithm is perfectly suitable for distributed storage and multi-core machines. In particular, since the number of reads $N$ is much larger than the number of sequences $M$, one can store the reads $R_1,\ldots,R_N$ in $n_p$ different processing units. Hence at each iteration $r$, the processing unit $p, p=1,\ldots,n_p$, can compute the local value
\[
\hat{\rho}_{m,p}^{r+1} = \frac{1}{N} \sum_{n \in \mathcal{N}_p} \frac{ \alpha_{nm} \rho_m^r}{ \sum_{m^\prime=1}^M  \alpha_{nm^\prime} \rho_{m^\prime}^r}, \quad  \forall m=1,\ldots, M,
\]
where $\mathcal{N}_p$ is the set of reads stored at processor $p$ with $\cup_{p=1}^{n_p} \mathcal{N}_p = \{1,2,\ldots,N\}$. Then, all processors  update their global abundance estimate through the consensus procedure
\[
\rho_m^{r+1} = \sum_{p=1}^{n_p} \hat{\rho}_{m,p}^{r+1}, \quad  \forall m=1,\ldots, M.
\]
For a very recent application of BSUM algorithm in gene RNA-seq abundance estimation, the readers are referred to \cite{bray2015near}.
 }

\subsubsection{Tensor decomposition}
The CANDECOMP/PARAFAC (CP) decomposition has applications in different areas such as 
 clustering \cite{shashua2005non} and compression \cite{ibraghimov2002application}. For ease of presentation, here we only consider order three tensors. Given a third order tensor $\mathfrak{X} \in \mathbb{R}^{m_1 \times m_2  \times m_3}$, its rank~$R$ CANDECOMP/PARAFAC decomposition is given by
$
\mathfrak{X} = \sum_{r=1}^R  \mathbf{a}_{r} \circ \mathbf{b}_{r}  \circ\mathbf{c}_{r},
$ where  $\mathbf{a}_{r} \in \mathbb{R}^{m_1}$, $\mathbf{b}_{r} \in \mathbb{R}^{m_2}$,  $\mathbf{c}_{r} \in \mathbb{R}^{m_3}$; and the notation ``$\circ$" stands for the outer product operator.

In general, finding the CP decomposition of a given tensor is NP-hard \cite{haastad1990tensor}. A well-known algorithm for finding the CP decomposition is the alternating least squares (ALS) algorithm proposed in \cite{carroll1970analysis,harshman1970foundations}.  This algorithm is in essence the BCD algorithm on the following optimization problem
\begin{equation}
\label{eq:ALSTensor}
\min_{\{\mathbf{a}_r,\mathbf{b}_r,\mathbf{c}_r\}_{r=1}^R} \quad \|\mathfrak{X} - \sum_{r=1}^R  \mathbf{a}_{r} \circ \mathbf{b}_{r}  \circ\mathbf{c}_{r}\|_F^2.
\end{equation}
In the ALS algorithm, we consider three blocks of variables: $\{\mathbf{a}_r\}_{r=1}^R$, $\{\mathbf{b}_r\}_{r=1}^R$, and $ \{\mathbf{c}_r\}_{r=1}^R$. At each iteration of the  algorithm, two blocks are held fixed and only one block is updated by solving~\eqref{eq:ALSTensor}. The block selection rule is cyclic and therefore one needs the uniqueness of the minimizer assumption at each iteration for theoretical convergence guarantee.  Clearly, this assumption does not hold in general in \eqref{eq:ALSTensor} and therefore, theoretically, convergence is not always guaranteed. In addition, another well-known drawback of the ALS algorithm is the {\it ``swamp"} effect where the objective remains almost constant for many iterations and then starts decreasing again. It has been observed in the literature that the employment of proximal upper-bound (see Table~\ref{table:bounds}) could help reducing the swamp effect \cite{CPSwampTikhonov}. It is also suggested in \cite{CPSwampTikhonov} that decreasing the proximal coefficient ($\gamma$ in Table~\ref{table:bounds}) during the ALS algorithm could further improve the performance of the  algorithm. Notice that these modifications in the algorithm makes the algorithm a special case of BSUM framework.  Consequently, its theoretical convergence is also guaranteed by Theorem~\ref{thm:main}.

Figure~\ref{fig:ALSTensor} compares the performance of the naive ALS algorithm with the one  using proximal upper-bound.  As can be seen from the figures the proximal ALS algorithm has less swamp effect as compared to naive ALS method. The  reader is referred to \cite{Razaviyayn12SUM,CPSwampTikhonov} for more details of the algorithm and to \cite{bazerque2013rank} for the application of BSUM and CP decomposition in gene expression and brain imaging.

 \begin{figure*}[ht]
 \begin{minipage}[t]{0.48\linewidth}
\includegraphics[width=1\linewidth]{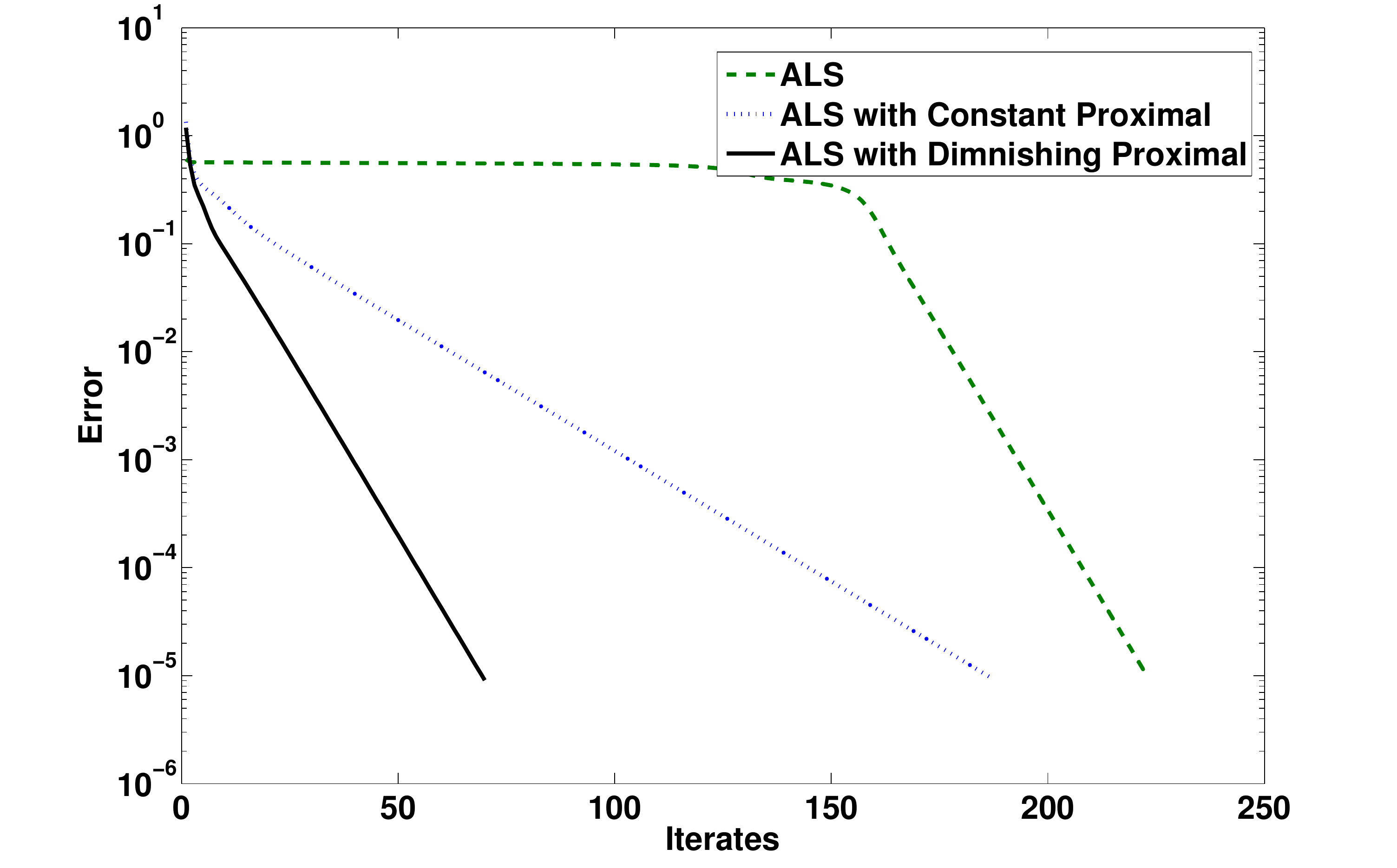}
\end{minipage}\hfill
 \begin{minipage}[t]{0.48\linewidth}
\includegraphics[width=1\linewidth]{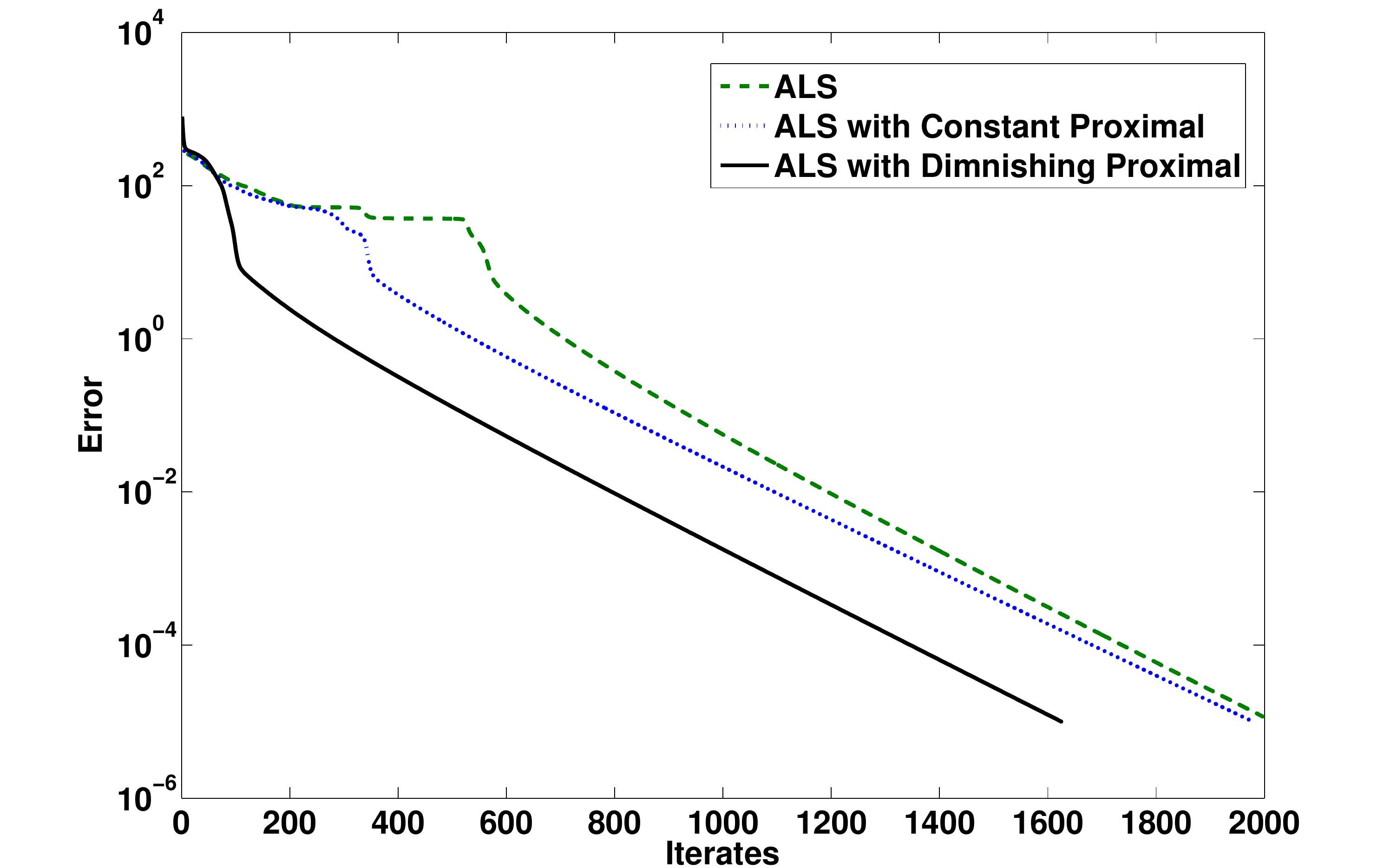}
\end{minipage}
\caption{\footnotesize \cite{Razaviyayn12SUM} Comparison of the ALS and Proximal ALS algorithm. The proximal ALS algorithm is the BSUM approach using the proximal upper bound - see Table~\ref{table:bounds}. The ``ALS with Diminishing Proximal" algorithm utilizes a decreasing proximal coefficient during the iterates of the algorithm. In the example $m_1=m_2=m_3=R=100$. \protect\footnotemark}\label{fig:ALSTensor}
\end{figure*}

\footnotetext{Copyright $\copyright$2013 Society for Industrial and Applied Mathematics.  Reprinted with permission.  All rights reserved.}

\subsection{Machine learning: sparse dictionary learning and sparse linear discriminant analysis}

\subsubsection{Dictionary learning for sparse representation}
{In compressive sensing  \cite{donoho2006compressed, candes2006robust} problems, a given data signal is represented by sparse linear combination of the signals in a given set called {\it dictionary}. In many applications even the dictionary is not known {\it a priori}, therefore it should be learned from the data.} More precisely, given a set of training signals $\{\mathbf{y}_1,\ldots,\mathbf{y}_N \in \mathbb{R}^n\}$, the dictionary learning task it to find a dictionary set $\{\mathbf{a}_1,\ldots, \mathbf{a}_k \in \mathbb{R}^n\}$ that can sparsely represent the signals in the training set. Defining the matrices $\mathbf{Y} \triangleq [\mathbf{y}_1,\ldots, \mathbf{y}_N]$, $\mathbf{A} \triangleq [\mathbf{a}_1 ,\ldots,\mathbf{a}_k]$, and $\mathbf{X} \triangleq [\mathbf{x}_1,\ldots, \mathbf{x}_N]$, the dictionary learning problem can be written as \cite{razaviyayn2014dictionary,lee2006efficient}
\begin{equation}
\nonumber
\begin{split}
\min_{\mathbf{A},\mathbf{X}} \quad &d(\mathbf{Y},\mathbf{A},\mathbf{X})\\
{\rm s.t.}\quad &\mathbf{A} \in \mathcal{A}, \;\mathbf{X} \in \mathcal{X},
\end{split}
\end{equation}
 where the sets $\mathcal{X}$ and $\mathcal{A}$ are given based on the prior knowledge on the data. The function $d(\cdot,\cdot,\cdot)$  measures the goodness-of-fit of the model. For example, a popular choice of the function $d(\cdot,\cdot,\cdot)$ and the set $\mathcal{A}$ leads to the following optimization problem  \cite{lee2006efficient}
 \begin{equation*}
\begin{split}
\min_{\mathbf{A},\mathbf{X}} \quad & \|\mathbf{Y} - \mathbf{A} \mathbf{X}\|^2 + \lambda \|\mathbf{X}\|_1\\
{\rm s.t.}\quad &\|\mathbf{a}_i\|_F^2 \leq \beta_i,
\end{split}
\end{equation*}
where the first term in the objective keeps our estimated signals close to the training set and the second term forces the representation to be  sparse.  One popular approach in the dictionary learning algorithm is to alternatingly update the dictionary $\mathbf{A}$ and the coefficients $\mathbf{X}$ \cite{aharon2006svd} . However, naively updating these two variables to its global optimum  requires solving a  sparse recovery problem at each iteration which is costly for big size problems. Motivated by the idea of inexact steps in BSUM framework, one can iteratively replace the objective by a locally tight upper-bound which is easier to minimize at each iteration and hence leading to computationally cheaper steps in the algorithm. It is not hard to see that utilizing the  quadratic upper-bound in Table~\ref{table:bounds} with diagonal matrices $\bfPhi_i$ leads to closed form updates at each step  \cite{razaviyayn2014dictionary}. Unlike many existing algorithms in the literature \cite{aharon2006svd,lee2006efficient}, the resulting algorithm is guaranteed to converge to the set of stationary solutions theoretically as the result of Theorem~\ref{thm:main}.

Figure~\ref{fig:DL} and Table~\ref{table_psnr} show the performance of the resulting algorithm for dictionary learning in an image denoising problem. The denoising is performed on the Lena image corrupted by additive Gaussian noise with various variances~$\sigma^2$. As can be seen from Table~\ref{table_psnr}, the proposed algorithm results in larger PSNR values than the K-SVD method \cite{aharon2006svd} when the noise level is large. Moreover, the proposed algorithm contains less visual artifacts. Furthermore,  each step of the proposed algorithm is in closed form and it is computationally favorable, while each step of the K-SVD method requires an inner iterative method.

\begin{figure}
  \centering
  \includegraphics[width=10.6cm]{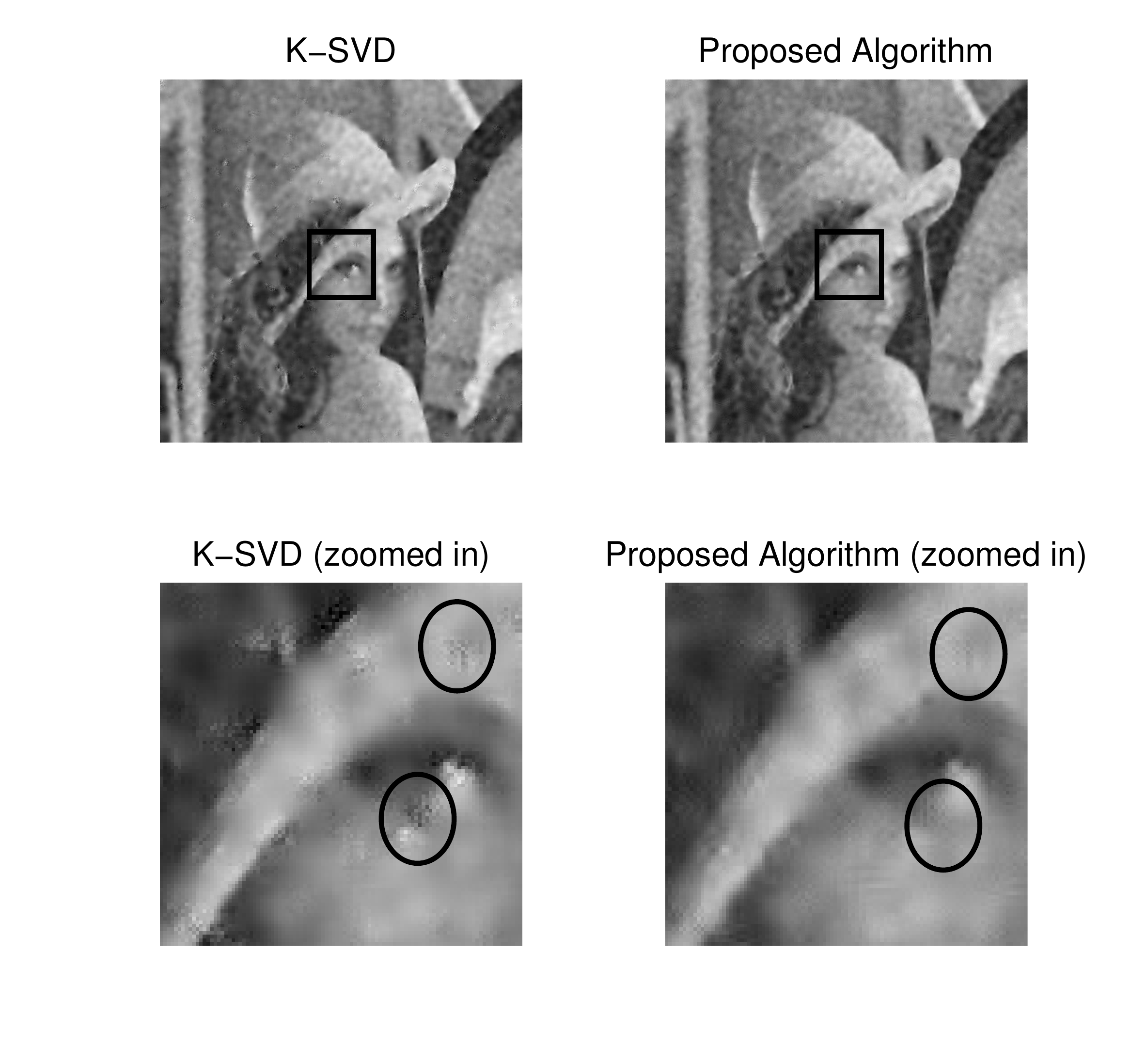}\\
  \vspace{-0.9cm}
  \caption{ \cite{razaviyayn2014dictionary} Sample denoised images ($\sigma = 100$). }\label{fig:DL}
  \vspace{-0.2cm}
\end{figure}
\begin{table}[h]
\centering
\small
\vspace{-0.3cm}
\begin{tabular}{|c|c|c|c|}
  \hline
  {\small$\sigma$/PSNR } & {\small  DCT } & {\small  K-SVD } & {\small  Proposed Algorithm} \\
  {\small 20/22.11} & {\small 32 } & {\small  \textbf{32.38}    } & {\small  30.88} \\
  {\small 60/12.57} & {\small  26.59 } & {\small  \textbf{26.86} } & {\small  26.37} \\
  {\small 100/8.132} & {\small  24.42 } & {\small  24.45 } & {\small  \textbf{24.46}} \\
  {\small 140/5.208} & {\small  22.96 } & {\small  22.93 } & {\small  \textbf{23.11} } \\
  {\small 180/3.025} & {\small  21.73 } & {\small  21.69 } & {\small  \textbf{21.96} }\\
  \hline
\end{tabular}
\caption{  \cite{razaviyayn2014dictionary} Image denoising result comparison on ``Lena image" for different noise levels. Values are averaged over $10$ Monte Carlo simulations.}
\label{table_psnr}
\end{table}

\subsubsection{Sparse linear discriminant analysis}
The linear discriminant analysis (LDA), which is closely related to analysis of variance (ANOVA) and regression analysis, is widely used in machine learning and statistics for classification and dimensionality reduction purposes; see, e.g., \cite{mclachlan2004discriminant}.  Let us, for the ease of presentation, focus only on the binary classification problem: Let $\mathbf{x}_i \in\mathbb{R}^p$, $i=1,2,\ldots,N$, denote the zero-centered observations, where each observation $\mathbf{x}_i$ belongs to one and only one of the two classes $\mathcal{C}_0$ and $\mathcal{C}_1$. Given the binary classes, the standard {\it within class covariance} estimate can be calculated by
\[
\widehat{\Sigma}_w = \frac{1}{N} \sum_{k\in\{0,1\}} \sum_{i \in \mathcal{C}_k} (\mathbf{x}_i - \widehat{\mu}_k) (\mathbf{x}_i - \widehat{\mu}_k)^T,
\]
where $\widehat{\mu}_k = \frac{1}{N} \sum_{i \in \mathcal{C}_k} \mathbf{x}_i$ is observations mean in class $\mathcal{C}_k$. Similarly, the standard {\it between class covariance} estimate is given by
\[
\widehat{\Sigma}_b = \frac{1}{N} \left(N_0 \widehat{\mu}_0 \widehat{\mu}_0^T + N_1 \widehat{\mu}_1 \widehat{\mu}_1^T\right),
\]
with $N_0$  (resp. $N_1$) being the cardinality of the set $\mathcal{C}_0$ (resp. $\mathcal{C}_1$). The goal of LDA is to find a lower dimensional subspace so that the projection of the observations onto the selected subspace leads to well-separated classes. In other words, the task is to project data points into a subspace with large between class variance relative to the within class variance. Let us for simplicity consider projection onto one dimensional  subspace defined by the vector $\boldsymbol{\beta} \in \mathbb{R}^p$; see \cite{witten2011penalized} {\color{black}for  details on projection to larger dimensional subspaces.} Then the inner product $\langle \boldsymbol{\beta} , \mathbf{x}\rangle$ is the projection of the observation $\mathbf{x}$ onto the selected subspace; and the within class variance of the projected data points is given by $\widehat{\sigma}_w = \boldsymbol{\beta}^T \widehat{\Sigma}_w \boldsymbol{\beta}$;
while the  between class variance can be written as $\widehat{\sigma}_b = \boldsymbol{\beta}^T \widehat{\Sigma}_b \boldsymbol{\beta}.$
Therefore, in the linear discriminant analysis problem, we are interested in solving:
\begin{equation}
\label{eq:LDA}
\begin{split}
\max_{\boldsymbol{\beta}}  \quad&\boldsymbol{\beta}^T \widehat{\Sigma}_b \boldsymbol{\beta} \\
{\rm s.t.} \quad & \boldsymbol{\beta}^T \widehat{\Sigma}_w \boldsymbol{\beta} \leq 1.
\end{split}
\end{equation}
Unfortunately, when the number of features is large relative to $N$, the matrix $\widehat{\Sigma}_w$ is rank deficient and therefore problem~\eqref{eq:LDA} is ill-posed. In order to resolve this issue and to have a small generalization error, \cite{witten2011penalized} suggests to regularize the optimization problem with a convex penalty function $P(\cdot)$; and solve
\begin{equation}
\label{eq:SLDA}
\begin{split}
\max_{\boldsymbol{\beta}}  \quad&\boldsymbol{\beta}^T \widehat{\Sigma}_b \boldsymbol{\beta}  - P(\boldsymbol{\beta})\\
{\rm s.t.} \quad & \boldsymbol{\beta}^T \widehat{\Sigma}_w \boldsymbol{\beta} \leq 1.
\end{split}
\end{equation}
Clearly, this optimization problem is non-convex. As suggested in \cite{witten2011penalized}, one can linearize the first part of the objective in \eqref{eq:SLDA} iteratively to obtain a tight upper-bound of the objective. It is not hard to see that the algorithm used in \cite{witten2011penalized} is BSUM with the {\it linear upper-bound} given in Table~\ref{table:bounds}.

\section{Extensions}\label{sec:extension}
In this section, we discuss  extensions and generalizations of the BSUM framework in various settings.
\subsection{Stochastic optimization}
Consider the following stochastic optimization problem
\begin{equation}
\label{eq:StochasticOpt}
\begin{split}
\min_x \quad &f(x) \triangleq \mathbb{E}_\xi \left[g(x,\xi)\right]\\
{\rm s.t.} \quad & x \in \mathcal{X},
\end{split}
\end{equation}
where $\mathcal{X}$ is a closed convex set and $\xi$ is a random variable modeling the uncertainty in our optimization problem. A standard classical approach for solving \eqref{eq:StochasticOpt} is the sample average approximation (SAA) method; see \cite{kim2011guide} and the references therein. At iteration $r$ of SAA method, given a new realization $\xi^r$ of the random variable $\xi$, SAA method generates a new iterate $x^r$ by solving a problem {with the sample average $\frac{1}{r} \sum_{i=1}^r g(x,\xi^i)$ as its objective,
where $\xi^1,\xi^2,\ldots, \xi^r$ are independent identically-distributed realizations of the random variable $\xi$.}

A major drawback of the SAA method is that each of its iteration can be computationally very expensive. The computational inefficiency arises from either the non-convexity of the objective, or not having closed form solutions at each iteration.

Motivated by the BSUM framework, authors of \cite{razaviyayn2013stochastic,mairal2013stochastic} suggest to use an inexact version of SAA method, in which a sequence of upper bounds of the objective are minimized. In particular, at each iteration $r$, the optimization variable is updated by
\begin{equation*}
\label{eq:SSUM}
\begin{split}
x^r \in \arg\min_x \quad &\frac{1}{r} \sum_{i=1}^r \widehat{g}(x,x^{i-1},\xi^i)\\
{\rm s.t.} \quad & x \in \mathcal{X},
\end{split}
\end{equation*}
where $\widehat{g}(\cdot,x^{i-1},\xi^i)$ is an upper-bound of the function $g(\cdot,\xi^i)$ around the point $x^{i-1}$. The approximation function $\widehat{g}(\cdot,x^{i-1},\xi^i)$ is assumed to be in the form of BSUM approximation. The resulting algorithm, named stochastic successive upper bound minimization (SSUM), is guaranteed to converge to the set of stationary solutions almost surely; see \cite{razaviyayn2013stochastic} for more details.
Further, it is shown to be capable of dealing with various practical problems in signal processing and machine learning. For example, as we will see shortly, the authors in \cite{razaviyayn2013stochasticWMMSE} apply SSUM framework to cope with uncertainties in channel estimation for wireless  beamformer design problem. As another example, the   online sparse dictionary leaning algorithm proposed in \cite{mairal2010online} is a special case of SSUM.
%
%

The stochastic optimization framework is well suited for many modern big data applications, especially when the entire data set is not available initially and the data points are made available over time. These problems can be considered as the above general stochastic optimization problem; see also \cite{Slavakis14}. In addition, many statistical model fitting problems, such as the simple classical regression problem, 
 can be cast as minimizing the following sum-cost function $\sum_{\ell=1}^L g(x,\xi_\ell)$. Typically, the number of data points $L$ is very large, making it difficult for batch processing. Therefore, it is desirable to implement algorithms working with only one (or a few) data point at each step. In these scenarios, the stochastic optimization framework is useful since the sum-cost minimization problem can be viewed as a stochastic optimization problem $\min \; \mathbb{E} \left[g(x,\xi)\right]$ with $\xi$ being uniformly drawn from the set $\{\xi_1,\ldots,\xi_L\}$.

{As an example of the SSUM method, consider the wireless transceiver design problem discussed in Section~\ref{sec:WMMSE} where the channel coefficients $\{\mathbf{H}_{ij}\}_{i,j}$ are not exactly known. In this scenario, we can consider the channel coefficients as random variables and solve the following stochastic optimization problem
\begin{equation}
\label{eq:ExpectedSumRateMaximization}
\begin{split}
\max_{\mathbf{u},\mathbf{v}} \;\; &\sum_{k=1}^K \mathbb{E}\bigg[R_k (\mathbf{u},\mathbf{v}) \bigg]\\
{\rm s.t.} \quad & \|\mathbf{v}_k\|^2 \leq P_k, \quad \forall k=1,2,\ldots,K,
\end{split}
\end{equation}
which is the stochastic counterpart of the optimization problem \eqref{eq:SumRateMaximization}. Utilizing the upper bound~\eqref{eq:WMMSEUpperBound} in SSUM algorithm leads to the stochastic WMMSE algorithm~\cite{razaviyayn2013stochasticWMMSE}. Figure~\ref{fig:SWMMSE} illustrates the numerical performance of the SSUM methods as compared with SAA. At each iteration of SAA procedure, one should solve a nonconvex optimization problem.  Two different methods are considered: the gradient descent method with random initialization and the WMMSE algorithm which is known to converge in few iterations for this problem.  As can be seen from the figure, the running time of the SAA algorithm is much longer than that of the SSUM.

\begin{figure*}[ht]
\hspace{-0.75cm}
    \begin{minipage}[t]{0.55\linewidth}
    \centering
     {\includegraphics[width=
1\linewidth]{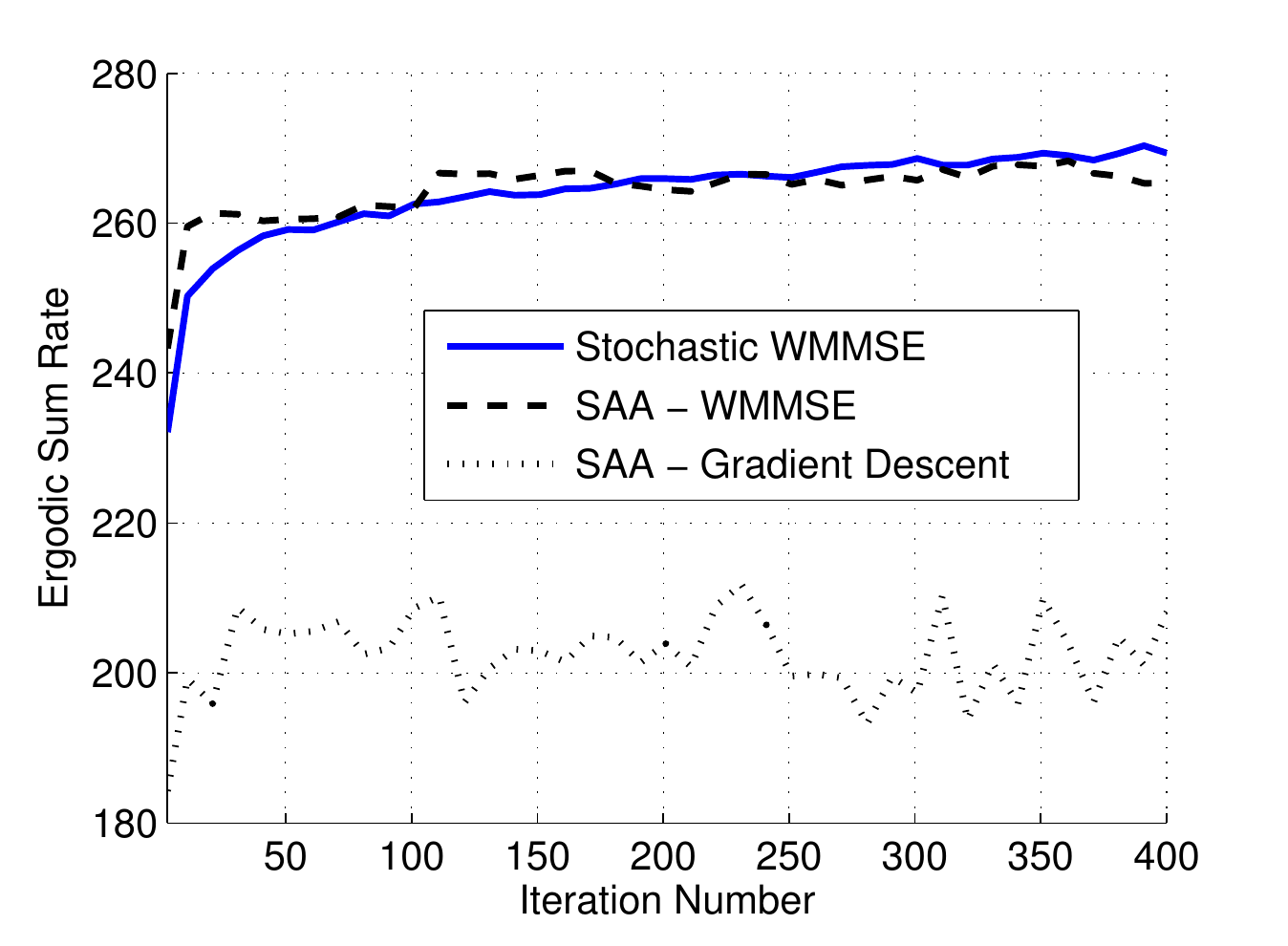}  }
\end{minipage}\hfill
    \begin{minipage}[t]{0.55\linewidth}
    \centering
    {\includegraphics[width=
1\linewidth]{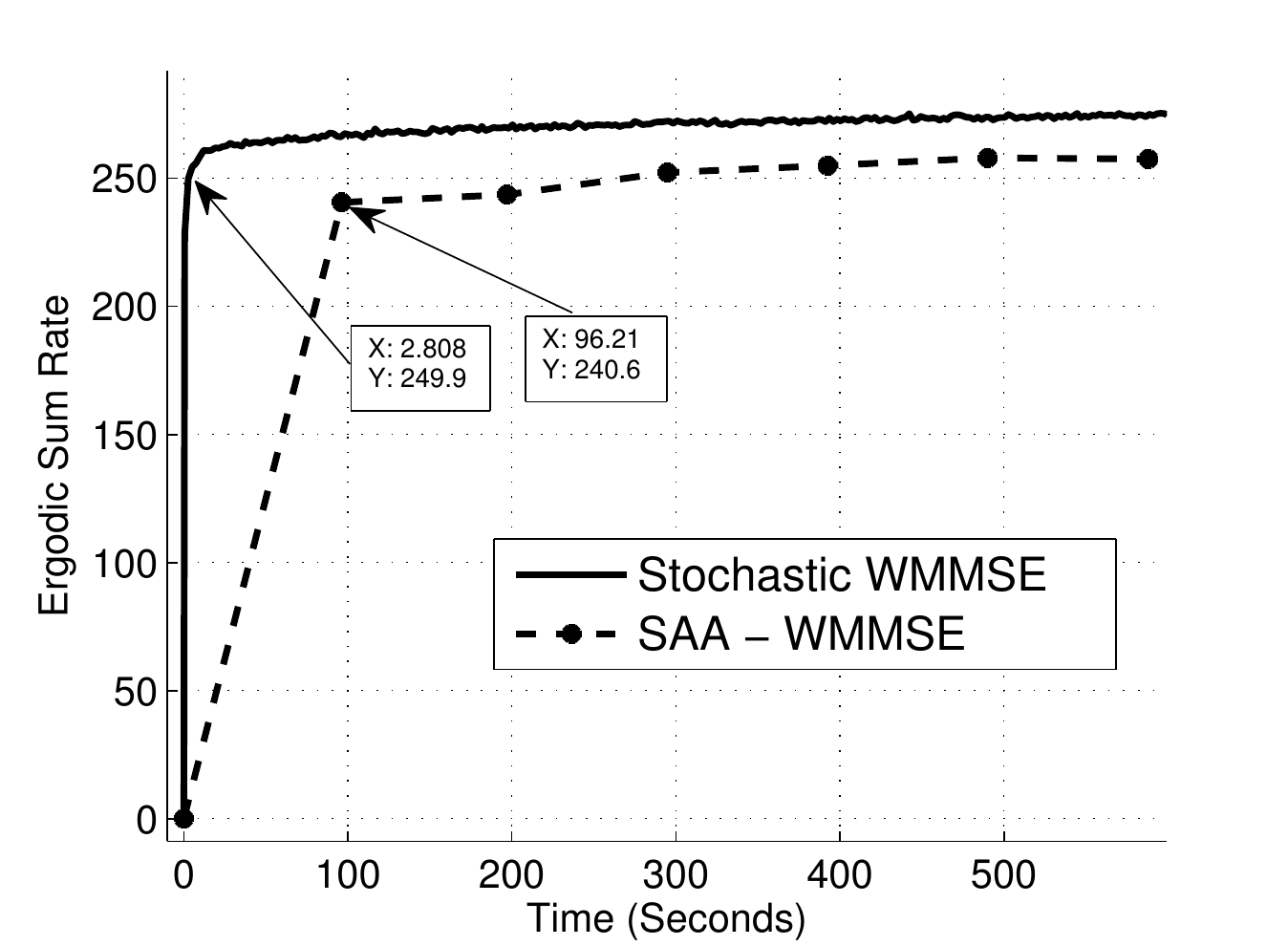}  }
\end{minipage}
\caption{Iteration and running time comparison of SSUM vs. SAA \cite{razaviyayn2013stochastic}}  \label{fig:SWMMSE}
    \end{figure*}

}

\subsection{Coupling constraints}
So far in this article, we have assumed that the constraints in the optimization problem is separable and convex. In other words, the constraint set $\mathcal{X}$ in \eqref{eq:OriginalProblem} is of the form $\mathcal{X} = \mathcal{X}_1 \times \ldots \mathcal{X}_n$ with each $\mathcal{X}_i$ being convex. A natural extension of the BSUM framework is to modify it in order to deal with coupling and nonconvex constraints. 
\subsubsection{Linear coupling}
Consider the following convex problem with linear coupling constraints
\begin{equation}
\label{eq:CouplingLinearProblem}
\begin{split}
\min_{(x_1,\ldots,x_n)} \; &f(x_1,\ldots,x_n) \\
\st \quad & \sum_{i=1}^n A_i x_i = b\\
& x_i \in \mathcal{X}_i, \;\;\forall i = 1,2,\ldots,n,
\end{split}
\end{equation}
where $x_i \in \mathbb{R}^{m_i}$, $A_i \in \mathbb{R}^{k \times m_i}$, $b \in \mathbb{R}^k$ and $f(\cdot)$ is a convex function. As seen in Example~\ref{ex:diverge:coupling}, the direct extension of the BCD/BSUM approach does not work for this type problem. A popular approach for solving the above optimization problem is the alternating direction method of multipliers (ADMM) \cite{BoydADMMsurvey2011,Glow84}.  This approach is based on finding a saddle point of the augmented Lagrangian function
\[
L(x_1,\ldots,x_n ; \lambda) = f(x_1,\ldots,x_n) + \left\langle  \lambda, \sum_{i=1}^n A_i x_i - b\right\rangle + \frac{\rho}{2} \left\|\sum_{i=1}^n A_i x_i - b\right\|^2,
\]
where $\lambda \in \mathbb{R}^k$ is the Lagrange multiplier corresponding to the linear constraint; $\rho>0$ is the augmented Lagrangian coefficient; and $\langle \cdot,\cdot \rangle$ denotes the inner product operator.

At each iteration of the ADMM method, either a primal block variable $x_i$ is  updated according to
\[
x_i^{r+1} \leftarrow \arg \min_{x_i \in \mathcal{X}_i} \quad L(x_i, x_{-i}^r;\lambda^r),
\]
or the dual Lagrange multiplier $\lambda$ is updated according to the gradient ascent rule
\[
\lambda^{r+1} \leftarrow \lambda^r + \alpha^r \left(\sum_{i=1}^n A_i x_i^r -b\right),
\]
where $\alpha^r$ is the dual step-size at iteration $r$. The  update orders for the primal and dual variables could be either cyclic or randomized.

Similar to the BSUM framework, 
one can replace the augmented Lagrangian function $L(\cdot,x_{-i}^r; \lambda^r)$ with its tight upper-bound $\widetilde{L}_i (\cdot,x^r;\lambda^r)$ at iteration $r+1$ where
\[
\widetilde{L}_i (x_i,x^r,\lambda^r) = u_i(x_i,x^r) + \left\langle \lambda^r,A_ix_i + \sum_{j\neq i} A_j x_j^r - b\right\rangle + \frac{\rho}{2} \left\|A_ix_i + \sum_{j\neq i} A_j x_j^r - b\right\|^2
\]
with $u_i(\cdot,x^r)$ being a locally tight approximation of the function $f(\cdot,x_{-i}^r)$ around the point $x_i^r$ satisfying Assumption A. The resulting algorithm, named Block Successive Upper Bound Minimization
Method of Multipliers  (BSUMM) \cite{hong13BSUMM},  is guaranteed to converge to the global optimal of \eqref{eq:CouplingLinearProblem} under some regularity assumptions \cite{hong13BSUMM}.
For extensions to non-convex problems, please see the recent work \cite{hong14nonconvex_admm,Ames13LDA}. {There are a few other interesting techniques to deal with linearly coupling constraint. For example references \cite{Necoara11, Reddi15} propose to randomly pick {2} blocks of variables to update at each iteration, and reference \cite{aybat15} propose new algorithms based on minimizing the augmented Lagrangian function. We refer the readers to these papers for more related works in this direction.  }

{
Let us illustrate an application of BSUMM to a multi-commodity routing problem, which arises in the design of next generation cloud-based communication networks \cite{liao15semi-async}. Consider a connected wireline network $\mathcal{N}=(\mathcal{V},\mathcal{L})$ which is controlled by $K+1$ Network Controllers (NCs) as illustrated in Fig. \ref{NetworkIllustration}. Let $\mathcal{V}$ denote the set of network nodes, which is partitioned into $K$ subsets, i.e., $\mathcal{V}=\cup_{i=1}^{K}\mathcal{V}^{i}$, $\mathcal{V}^{i}\cap\mathcal{V}^{j}=\emptyset$, $\forall\;i\neq j$. The set of directed links that connect nodes of $\mathcal{V}$ is denoted as $\mathcal{L}\triangleq\{l=(s_{l},d_{l})\mid\;\forall\;
s_{l},d_{l}\in\mathcal{V}\}$, where $l=(s_{l},d_{l})$ denotes the directed link from node $s_{l}$ to node $d_{l}$. Each NC $i$ controls $\mathcal{V}^{i}$ and the links connecting these nodes, i.e., $\mathcal{L}^{i}\triangleq\{l=(s_{l},d_{l})\in\mathcal{L}\mid \forall\;s_{l},d_{l}\in\mathcal{V}^{i}\}$ (cf. Fig. \ref{NetworkIllustration}). The network $\mathcal{N}^{i}\triangleq(\mathcal{V}^{i},\mathcal{L}^{i})$ is called the subnetwork $i$.

Our objective is to transport $M$ data flows over the network, with each flow $m$ being routed from a source node $s_{m}\in\mathcal{V}$ to the sink node $d_{m}\in\mathcal{V}$. We use $r_{m}\geq 0$ to denote flow $m$'s rate, and use $f_{l,m}\geq 0$ to denote its rate on link $l\in\mathcal{L}$. We also assume that a master node exists which controls the data flow rates $\{r_{m}\}_{m=1}^{M}$.  The central NC $0$ controls the subnetwork $\mathcal{N}^{0}$, consisting of the master node and the links connecting different subnetworks, i.e., $\mathcal{L}^{0}=\cup_{i\neq j}\mathcal{L}_{ij}^{0}$.

Let us consider two types of network constraints, as we list below.
\begin{enumerate}

\item {\it Link capacity constraints}. Assume each link $l\in\mathcal{L}$ has a fixed capacity denoted as
$C_{l}$. The total flow rate on link $l$ is constrained by\vspace{-0.1cm}
\begin{align}\label{CapacityWired}
{\bf 1}^{T}\mathbf{f}_{l}\leq C_{l},~\forall\;l\in\mathcal{L},
\end{align}
where ${\bf 1}$ is the all-one vector and $\mathbf{f}_{l}\triangleq[f_{l,1},\ldots,f_{l,M}]^{T}$.

\item {\it Flow conservation constraints}.  For any node $v\in\mathcal{V}$ and data flow $m$, the total incoming flow should be equal to the total outgoing flow:\vspace{-0.1cm}
\begin{align}\label{ConservationConst}
&\sum_{l\in{\rm In}(v)}f_{l,m}+1_{v=s(m)}r_{m}=\sum_{l\in{\rm Out}(v)}f_{l,m}+1_{v=d(m)}r_{m},\quad~\; m=1,\cdots,
M,~\forall\; v\in\mathcal{V},
\end{align}
where ${\rm In}(v)\triangleq\{l\in\mathcal{L}\mid\;d_{l}=v\}$ and ${\rm Out}(v)\triangleq\{l\in\mathcal{L}\mid\;s_{l}=v\}$ denote the set of
links going into and coming out of a node $v$ respectively; $1_{v=x}=1$ if $v=x$, otherwise $1_{v=x}=0$.
\end{enumerate}

To provide fairness to the users, let us maximize the minimum rate of all data flows.
The problem can be formulated as the following linear program (LP)\vspace{-0.3cm}
\begin{subequations}\label{MainProb}
\begin{align}
\max_{\mathbf{f},~\mathbf{r}}\quad r_{\min}~~~~{\rm s. t.}&\quad \mathbf{f}\geq 0,~r_{m}\geq r_{\min},~m=1,\cdots, M\label{NonNegativeConst}\\
&\quad\eqref{CapacityWired}~\mbox{and}~\eqref{ConservationConst},
\end{align}
\end{subequations}
where $\mathbf{f}\triangleq\{\mathbf{f}_{l}\mid l\in\mathcal{L}\}$ and $\mathbf{r}\triangleq\{r_{\min},r_{m}\mid m=1\sim M\}$. Obviously, one can use off-the-shelf optimization packages such as Gurobi  \cite{gurobi} to solve the LP \eqref{MainProb}, but this is only viable in a centralized setting where all the flows are managed by a single controller.


To enable distributed/parallel network management across the NCs, we need to allow each NC $i$ to independently optimize the variables belonging to the subnetwork $\mathcal{N}^i$. However, this task is difficult because the optimization variables of problem \eqref{MainProb} is {\it coupled} (indeed each flow rate $f_m$ appears in exactly {\it two} flow conservation constraints). To address this problem, we introduce a few sets of new variables to decouple the flow conservation constraints across different subnetworks (we refer the readers to \cite{liao15semi-async} for the detailed reformulation). The reformulated problem $\eqref{MainProb}$ is given by\vspace{-0.3cm}
\begin{align*}
\max_{\bx}&\quad r_{\min}\nonumber\\
{\rm s. t.}&\; {\{r_{\min},\bx_{02}\}\in\mathcal{X}_{0}},~{\{\bx_{i1},\bx_{i2}\}\in\mathcal{X}_{i1},~\bx_{i3}\in\mathcal{X}_{i2}},\\
&\underbrace{{\bx_{01}^{i}}={\bx_{02}^{i}}}_{\footnotesize\mbox{in $\mathcal{N}^{0}$}},~\underbrace{{\bx_{i1}}={\bx_{01}^{i}}}_{\footnotesize\mbox{in $\mathcal{N}^{i}$ and $\mathcal{N}^{0}$}},~\underbrace{{\bx_{i2}}={\bx_{i3}}}_{\footnotesize\mbox{in $\mathcal{N}^{i}$}},~i=1\cdots K\nonumber
\end{align*}
where $\mathcal{X}_0$ and $\mathcal{X}_{i1}$ and $\mathcal{X}_{i2}$ are some feasible sets, and $\{r_{\min},\bx_{02}, \bx_{01},\bx_{i1},\bx_{i2},\bx_{i3}\}$ are the block-variables. By applying the BSUMM, we can obtain a parallel/distributed algorithm. A few remarks about the implementation of this algorithm are in order.
\begin{enumerate}
\item The replication of link/flow variables for links across different subnetworks allows each subnetwork to be considered separately and independently. This feature makes the BSUMM subproblems solvable in parallel. The requirement of the replicated variables being the same as the original variables is enforced by the linear coupling constraints, and they can be satisfied asymptotically as the BSUMM algorithm converges.
\item The subproblems of the proposed BSUMM-based algorithm can be solved by each NC very efficiently. For example the update of $\{{r_{\min}},\bx_{02}\}$ can be performed by each NC in closed-form; the update of $\{{\bx_{i}},\bx_{01}^{i}\}$ can be performed by running the well-known RELAX code \cite{Bertsekas87}.
\item A careful implementation of the BSUMM allows the NCs to act asynchronously, in the sense that they do not need to coordinate with each other for computation. Such asynchronous implementation has the potential of greatly improving the computational efficiency.
\end{enumerate}

We illustrate the BSUMM implementation over a mesh wireline network with $126$ nodes, which is randomly partitioned into $9$ subnetworks with $306$ directed links within these subnetworks and $100$ directed links connecting the subnetworks. The capacities for the links within (resp. between) the subnetworks are uniformly distributed in [50,100] MBits/s (resp. [20,50] MBits/s). All simulation results are averaged over $200$ randomly selected data flow pairs and link capacity.

To demonstrate the benefit of parallelization, we also utilize a high performance computing (HPC) cluster, and make each computing node to be a NC. We compare a few different algorithms, listed below:
\begin{enumerate}
\item Solving a large-scale LP by Gurobi \cite{gurobi}, a centralized solver;
\item The synchronous BSUMM algorithm with $K=10$ NCs, with the computation done by either a single or by $10$ distributed computing cores;
\item The asynchronous BSUMM with $K=10$ NCs, {with the computation done in $10$ distributed computing cores.}
\end{enumerate}
Note that the asynchrony in the network arises naturally from the per-node computational delay and network communication delay. In the following table we demonstrate the performance of various algorithms when $M=200$. Clearly the asynchronous BSUMM with a small number of NCs outperforms all the rest of the algorithms.

\begin{table}[!htbp]
\centering\small
\caption{\normalsize Comparison with Different Algorithms. $M=200$}
\begin{tabular}{|c|c|c|}
\hline
{\bf Approaches} & {\bf Time} & {\bf \# of iterations}\\
\hline
Gurobi & $11.4690$s & N/A\\
Synchronous BSUMM ($10$ NCs, $1$ core) & $18.0679$s & 264.71\\
Synchronous BSUMM ($10$ NCs, $10$ core) & $5.80$s &  264.71\\
Asynchronous BSUMM ($10$ NCs, $10$ core) & $4.23$s & 2109.07\\
\hline
\end{tabular}
\end{table}
The numerical results suggest that appropriate network decomposition and asynchronous implementation are both critical for the fast convergence of BSUMM. In practice, we observe that the network should be decomposed following a few guidelines: {\it i)} the  computation burden across the subnetworks is well balanced; {\it ii)} the subroutine within the network can achieve its maximum efficiency; {\it iii)} the total number of replicated auxiliary variables is small.


\begin{center}
\vspace{-0.5cm}
\begin{figure}
    \centering
     {\includegraphics[width=
0.5\linewidth]{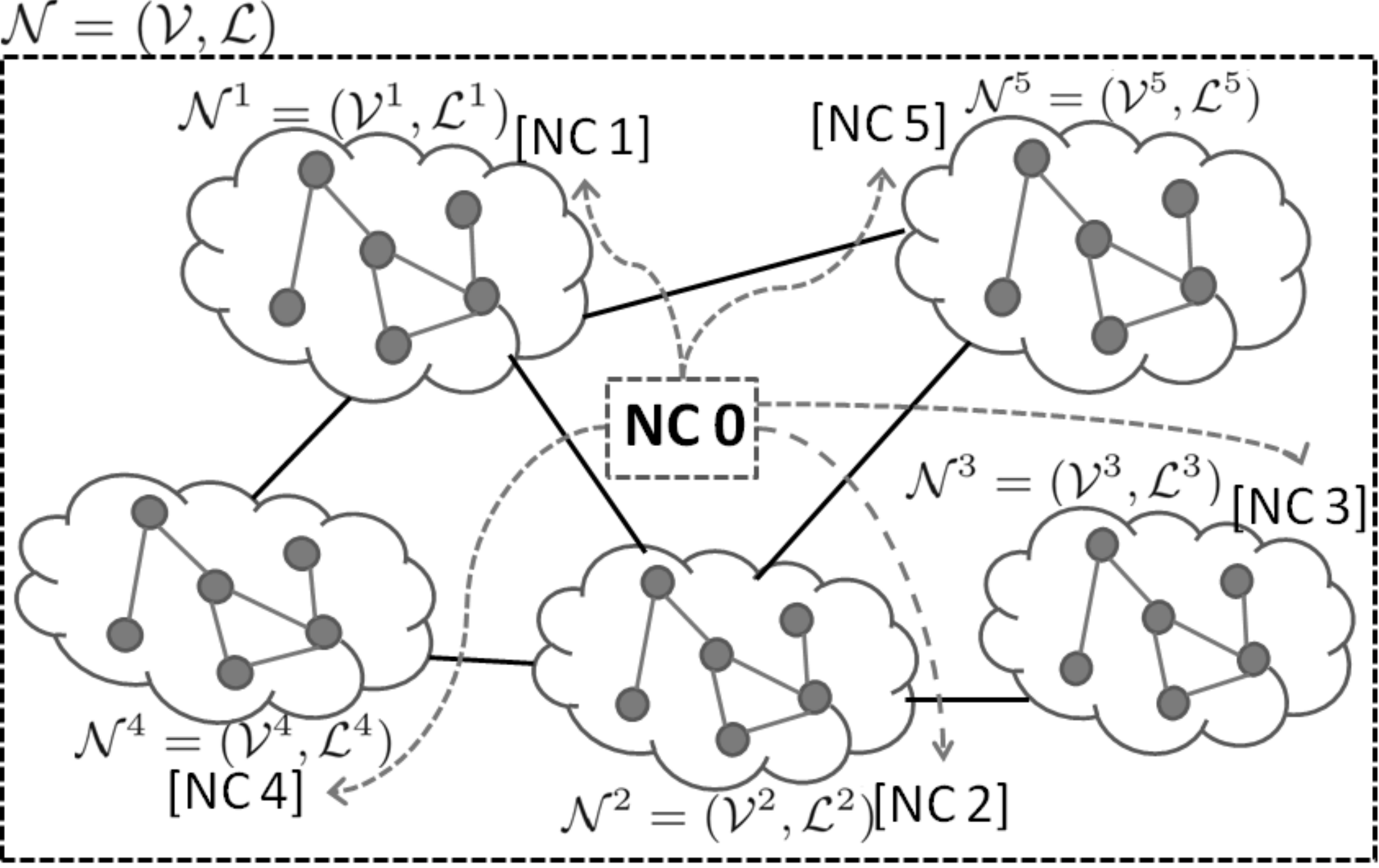}}\vspace{-0.3cm}
\caption{A wireline network consists of 5 subnetworks. Each of them is controlled by a NC, and these NCs are coordinated globally by a central NC 0 \cite{liao15semi-async}.}\label{NetworkIllustration}\vspace{-0.65cm}
\end{figure}
\end{center}
}
\subsubsection{Nonconvex constraints}
The  BSUM idea can be straightforwardly extended to a nonconvex constraint scenario for single block optimization problems. To proceed, let us consider the optimization problem
\begin{equation}
\label{eq:NonconvexConstProblem}
\begin{split}
\min_{x}  \quad&f_0(x) \\
\st \quad &  f_i(x) \leq 0,\quad \forall i = 1,2,\ldots,\ell,
\end{split}
\end{equation}
where the functions $f_i(\cdot)$ are not necessarily convex. Since dealing with nonconvex constraints is often not easy, one popular approach is to replace the functions $f_i(x), i = 1,2,\ldots,\ell,$ with their local tight upper-bound $u_i(x,x^r)$ iteratively. In other words, the update rule of the iterates is given by  \cite{Marks78}
\begin{equation}
\begin{split}
x^{r+1} \;\leftarrow\; \arg\min_{x} \quad &u_0(x,x^r) \\
\st \quad &  u_i(x,x^r) \leq 0,\quad \forall i = 1,2,\ldots,\ell.
\end{split}
\end{equation}
 As illustrated in Figure~\ref{fig:SCAconst}, the iterative approximation of the constraints is a restriction of the constraints and hence the iterates remain feasible during the algorithm. If, in addition, some constraint qualification conditions are satisfied, the resulting algorithm is guaranteed to converge to the set of stationary solutions of \eqref{eq:NonconvexConstProblem}; see \cite[Theorem 1]{razaviyayn14thesis} for detailed conditions and analysis.
 \begin{figure}[ht]
 \centering
\includegraphics[width=0.4\linewidth]{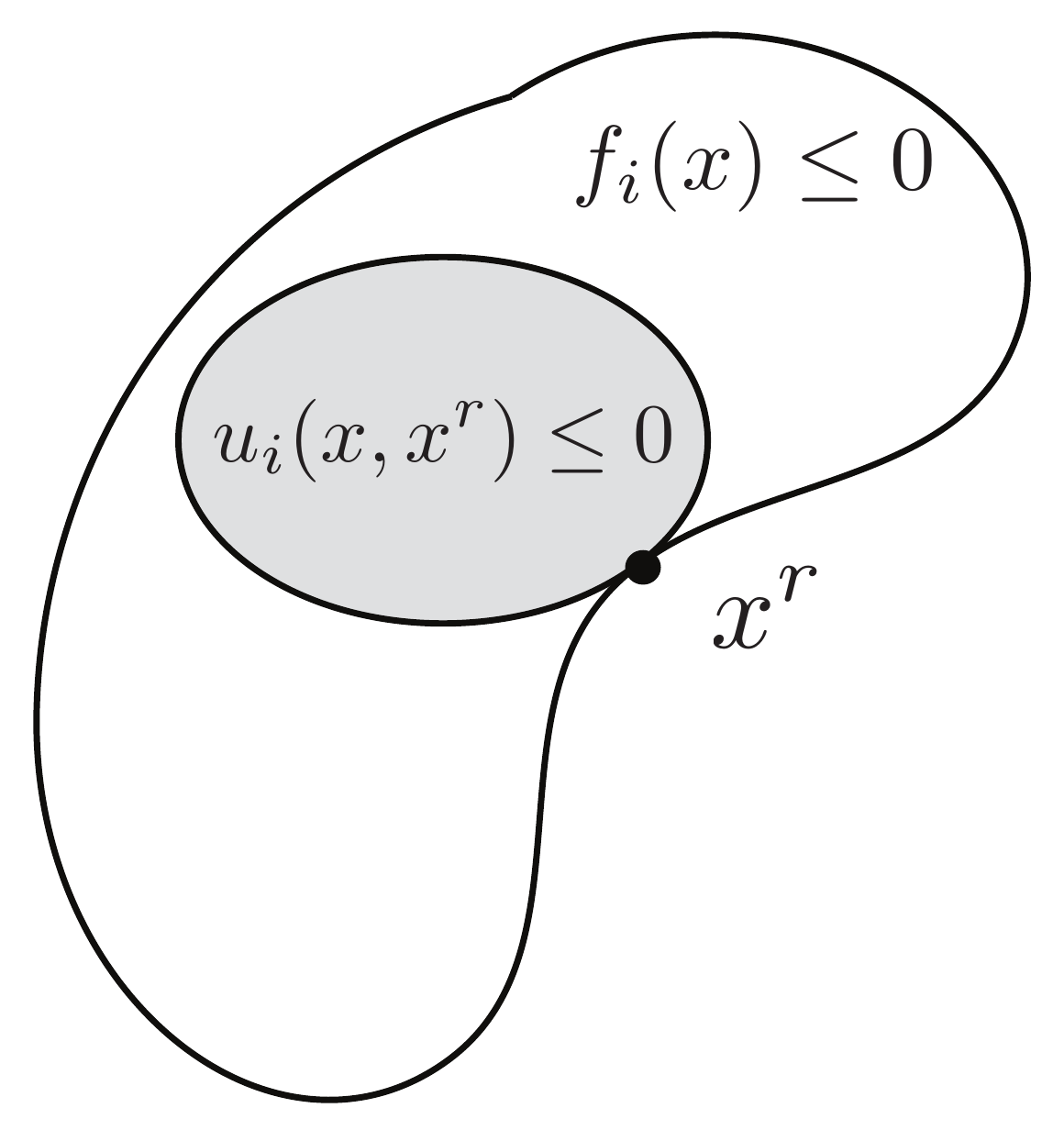}
\caption{Illustration of constraint convexification}
 \label{fig:SCAconst}
\end{figure}

\subsection{Parallel version and extensions to game theory}
\label{subsec:Parallel}
With the recent advances in multicore and cluster computational platforms, it is desirable to design {\it ``parallel"} algorithms for multi-block optimization problem where multiple cores update the block variables in parallel to optimize the objective function. A naive parallel extension of the BCD approach for solving \eqref{eq:OriginalProblem}  is to update all blocks (or a subset of them) in parallel by solving
\begin{equation}
\nonumber
x_i^{r+1} \;\leftarrow\; \arg \min_{x_i \in \mathcal{X}_i} \; \;\; f(x_i,x_{-i}^r), \quad \forall i = 1,\ldots,n.
\end{equation}
Unfortunately this naive extension of the BCD algorithm does not converge in general and might result in zig-zag/oscillating  or divergent behavior. As an example, consider the problem
\begin{equation*}
\begin{split}
\min_{(x_1,x_2)} \;\; \; &(x_1 - x_2)^2\\
\st \quad &-1 \leq x_1,x_2 \leq 1.
\end{split}
\end{equation*}
Clearly, this problem is convex with bounded feasible set and its optimal value is zero. However,  the above naive parallel extension of the algorithm leads to the following iteration path:
\[
(x_1^0,x_2^0) = (1,-1) \rightarrow (x_1^1,x_2^1) = (-1,1)\rightarrow (x_1^2,x_2^2) = (1,-1) \rightarrow \cdots,
\]
which is clearly not convergent. This is caused by aggressive steps used in the algorithm. To make the algorithm convergent, it is then necessary to employ small enough and controlled steps. Furthermore,  in the case of nonconvex objective function $f(\cdot)$ in \eqref{eq:OriginalProblem}, the approximation functions could be again used to obtain computationally efficient update rules. The resulting algorithm, dubbed parallel successive convex approximation (PSCA),  is summarized in Table~\ref{Table:ParallelBSUMAlgorithm}. To see the convergence  analysis of this algorithm and other related algorithms such as the Flexible Parallel Algorithm (FPA), the  reader is referred to \cite{meisam14nips, facchinei2013flexible, scutari13decomposition} and the references therein.

\begin{table}[htb]
\centering
\begin{tabular}{|p{4in}|}
\hline
\begin{itemize}
\item [1] \;Find a feasible point $x^0\in \mathcal{X}$; set $r = 0$, and choose a stepsize sequence $\{\gamma^r\}$
\item [2] \; \textbf{repeat}
\item [3] \quad {Pick index set $\cI^{\, r}$}
\item [4] \quad  {Let $\hat{x}_i^{r}= \arg \min_{x_i\in \mathcal{X}_i} u_i(x_i,x^{r-1})$, $\forall~i\in\cI^{\,r}$}
\item [5] \quad Set $x_k^{r} = x_k^{r-1}, \quad \forall\; k \notin \cI^{\,r}$
\item [6] \quad Set $x_i^{r} = x_i^{r-1} + \gamma^r (\hat{x}_i^r - x_i^{r-1}), \quad \forall\; i \in \cI^{\,r}$
\item [7] \quad $r = r+1$,
\item [8] \; \textbf{until} some convergence criterion is met
\end{itemize}
\\
\hline
\end{tabular}\vspace{1.2em}
\caption{Pseudo code of the PSCA algorithm}
\label{Table:ParallelBSUMAlgorithm} \vspace{-0.5cm}
\end{table}
Notice that PSCA can be viewed as a way of solving a multi-agent optimization problem where  multiple agents/users try to optimize a common objective by updating their own variable iteratively. Furthermore, it can  be  used in a game theoretic setting where each player in the game utilizes the best response strategy by optimizing a locally tight upper-bound of its own utility function. This algorithm is guaranteed to converge for some particular class of games  under some regularity assumptions on the players utility functions. The convergence analysis presented in  \cite{JSPMRGameTheory14} is based on certain contraction approach as well as monotone convergence for potential games.

Figure~\ref{fig:PSCA}  illustrates the behavior of cyclic and randomized parallel PSCA method as compared with their serial counterparts (i.e., the ``Cyclic BCD" and the ``Randomized BCD" ) applied to the LASSO problem. The performance of the ``PSCA" method is also illustrated for different number of processors and various block selection rules. As can be seen from the figures, parallelizaiton can result in more efficient algorithm; however, the convergence speed does not grow linearly with the number of processing cores. Moreover, increasing the number of processors beyond certain point results in slower convergence, which can be attributed to the increased communication overhead among the nodes.
\begin{figure*}[ht]
\hspace{-0.75cm}
    \begin{minipage}[t]{0.55\linewidth}
    \centering
     {\includegraphics[width=
1\linewidth]{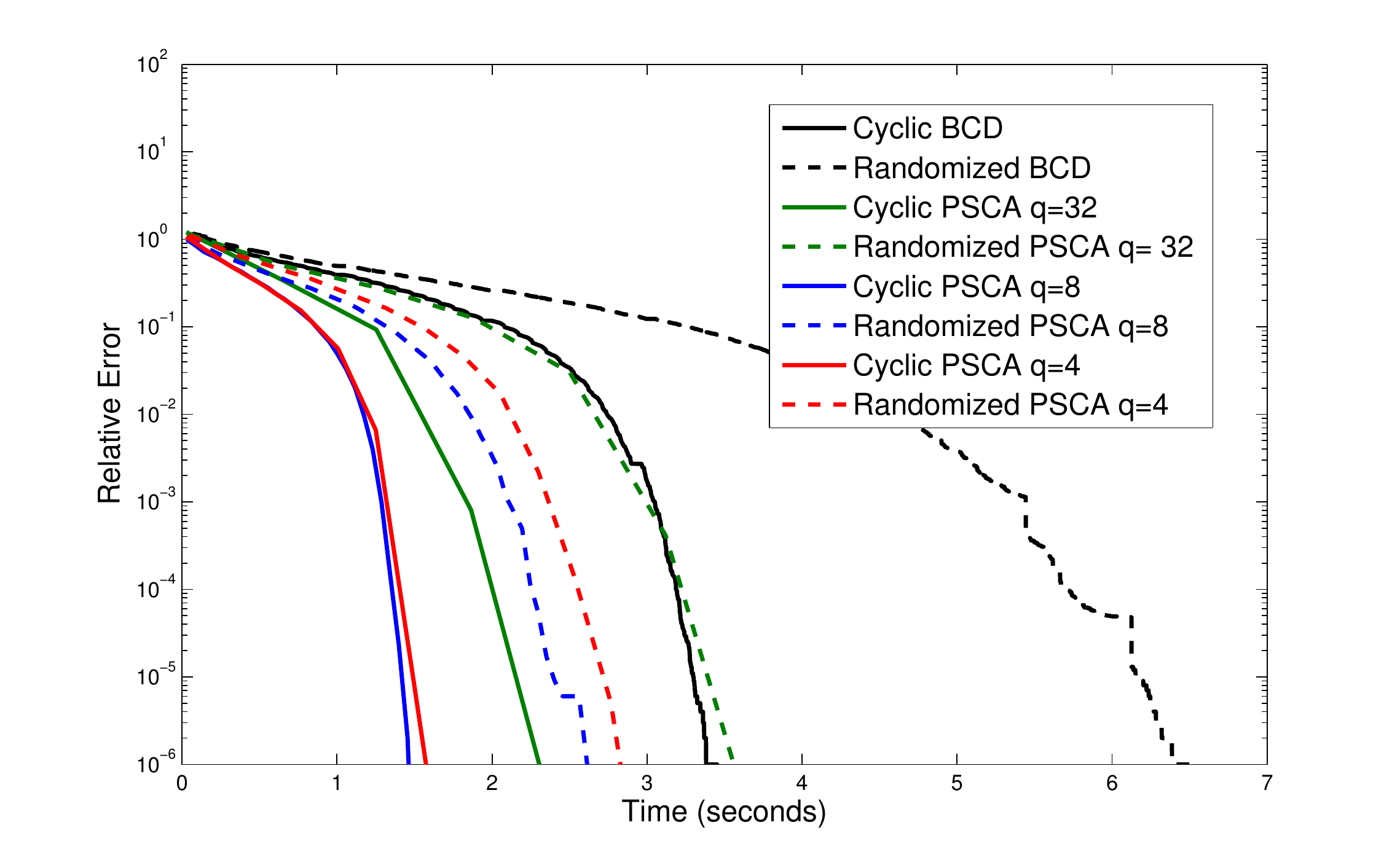}  }
\end{minipage}\hfill
    \begin{minipage}[t]{0.55\linewidth}
    \centering
    {\includegraphics[width=
1\linewidth]{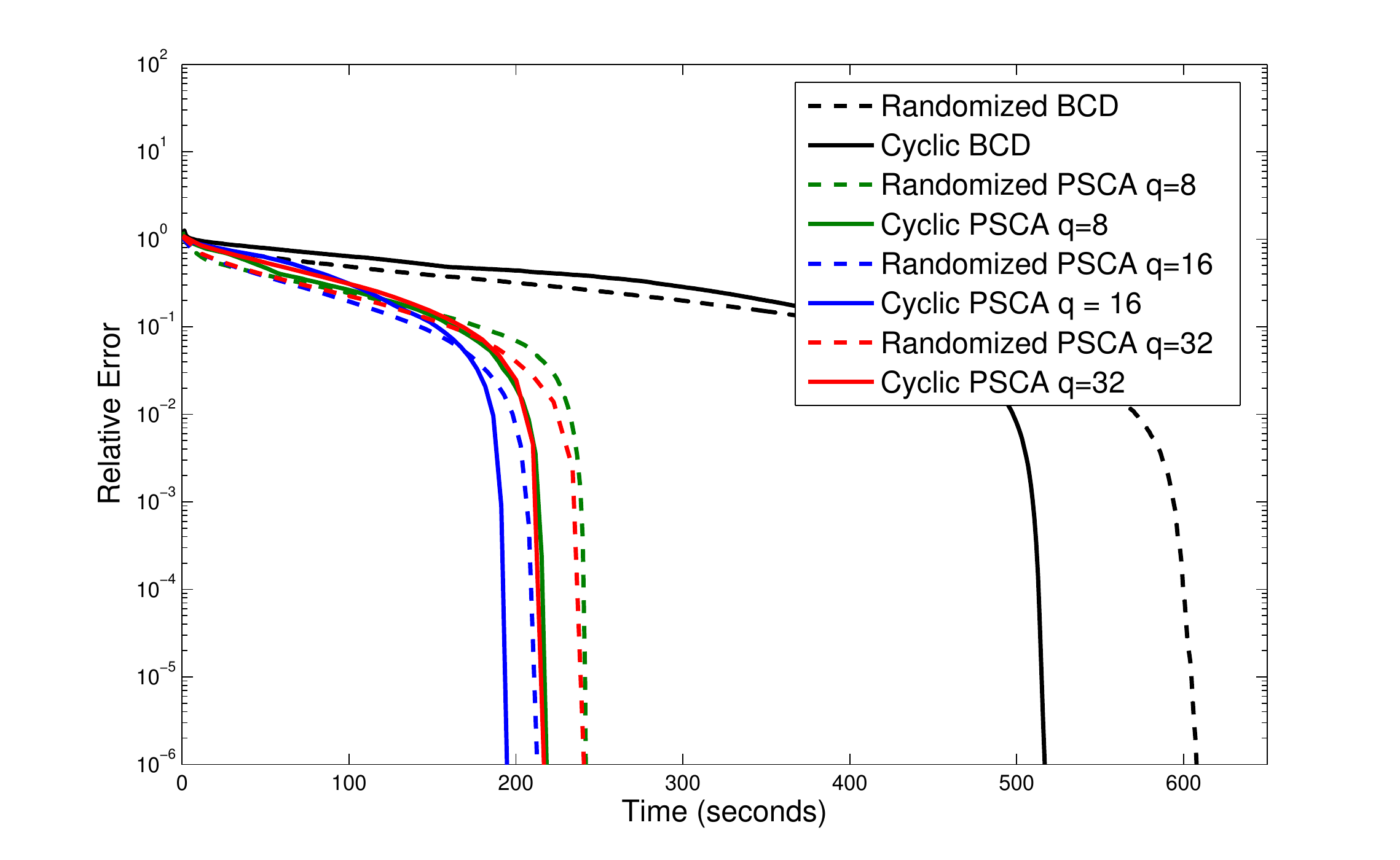}  }
\end{minipage}
\caption{ \cite{meisam14nips} Comparison of the serial BCD with PSCA method for solving the Lasso problem: $\min_x \; \|Ax - b\|^2 + \lambda \|x\|_1$. The matrix $A$ and vector $b$ are generated according to \cite{nesterov2013gradient}; and $q$ denotes the number of processors used in each experiment. The dimension of $A$ is  $2,000\times 10,000$ for the left figure and $1,000 \times 100,000$ for the right figure. {The experiments are done over a computer cluster using Message Passing Interface (MPI) in C.} }  \label{fig:PSCA}
    \end{figure*}

It is also worth noting that the parallel update rule is very useful in dealing with distributed data sets. Consider solving the LASSO problem
with the following objective: $\|Ax - b\|^2 + \lambda \|x\|_1.$ Assume we have $q$ processing cores each having their own memory.  Let us partition the matrix $A$ and vector $x$ into $q$ blocks: $A = [A_1,\ldots,A_q] \quad {\rm and }\quad x = [x^T_1,\cdots, x^T_q]^T$.
If PSCA is implemented in a way that each core $j$ is only responsible for updating block $j$, then at each iteration $r+1$, core~$j$'s problem of interest can be written as
\[
\min  \quad \|A_j x_j - b^\prime_j\|^2 + \lambda \|x_j\|_1,
\]
where $ b^\prime_j \triangleq b - \sum_{i \neq j} A_i x_i^r$. Notice that the value of $b^\prime_j$ can be calculated by letting each node $i$ compute the value of $A_i x_i$ and broadcast it to other nodes.  Hence under this architecture, each node does not need to know the complete matrix $A$ and only local information is enough for a distributed implementation of the PSCA method.

\subsection{\color{black}Practical Considerations}
There are a number of factors that we need to consider when using the BSUM framework.

The first consideration is about the choice of the upper bound. What is a good bound for a given application? The answer is generally problem dependent, as we have already seen in a few examples. The general guideline is that a good upper bound should be able to ensure algorithm convergence, best exploit the problem structure and make the subproblems easily solvable (preferably in closed-form). For example a simple proximal upper bound is not likely to perform well for the transceiver design problem discussed in Section \ref{sec:WMMSE}, as the resulting subproblems will not decompose over the variables.

The second consideration is about the choice of the block update rules. As we have seen in Section \ref{sub:BSUM:converge:rate}, different update rules lead to quite distinct convergence behaviour. For convex problems, deterministic rules such as the cyclic rule promise the worst-case convergence rates, while the randomized rule ensures convergence rate in either averaged or high probability sense.  Further, there is barely any theoretical rate analysis for nonconvex problems, regardless of the block selection rules. Therefore the best strategy in practice is to perform extensive numerical study and pick the best rule for the application at hand. For example, researchers have found that MBI rule is effective for certain tensor decomposition problem \cite{Chen2012MBI}; the cyclic rule can be superior to the randomized rule for certain LASSO problem, and certain G-So rule can outperform the randomized rule  \cite{Saha10,schmidt15icml}.

The third consideration is about the choice of the parallelization schemes. There has been extensive research on parallelizing various special cases and variations of the BSUM type algorithm, see \cite{razaviyayn2014parallel,richtarik13distributed,peng13,Scherrer12,liu14asynchronous, facchinei2013flexible, Facchinei15} and the references therein. These algorithms differ in a number of implementation details and in their applicability. For example, most of the implementations use randomized block selection rules to pick the variable blocks, while \cite{peng13} and \cite{facchinei2013flexible} additionally use certain variations of the G-So rule. The majority of the schemes only work for convex problems, with the exception of \cite{facchinei2013flexible,razaviyayn2014parallel} which work for general nonsmooth and nonconvex problems in the form of  \eqref{eq:problemBSUM2}. When assessing whether a given problem is suitable
for parallelization, it is important to know that oftentimes the number of blocks that can be updated in parallel is {\it data dependent}. For example when solving LASSO problems, it is shown in \cite{Scherrer12} and \cite{peng13} that the degree of parallelization is dependent on the maximum eigenvalue of certain submatrices of the data matrix. Some recent result \cite{Fercoq14} shows that for certain randomized coordinate descent method, such dependency could be mild. For solving general convex nonsmooth problems, \cite{richtarik13distributed} shows that the stepsize should be carefully selected based on both the ``separability" of the problem (or the sparsity of the data matrix) as well as the degree of parallelization. If the application at hand does not satisfy these conditions, the alternatives usually are: 1) exploit the problem structure and pick a good upper bound, so that the subproblems are decomposable, leading to parallel and {\it stepsize-free} updates; see for example the NMF problem in Section \ref{sub:nmf} and the WMMSE algorithm in Section \ref{sec:WMMSE}; 2) to use the diminishing stepsizes for updating the blocks, see for example \cite{facchinei2013flexible,razaviyayn2014parallel, Facchinei15}.

\section{Issues and Open Research Problems}
{This article presents a comprehensive algorithmic framework called BSUM for block-structured large-scale optimization. The main strength of the BSUM framework is its strong theoretical convergence guarantee and its flexibility. As demonstrated in this article, the BSUM framework covers a number of well-known but seemingly unrelated algorithms as well as their new extensions. Moreover, it is amenable to a number of different data models as well as to parallel implementation on modern multi-core computing platforms.

To close, we briefly highlight a couple of issues and open research topics related to the BSUM framework.
}
\begin{itemize}
\item {\bf Communication delay and overhead in parallel implementations}:
As discussed in  subsection~\ref{subsec:Parallel}, the convergence speed of the parallel version of BSUM framework does not increase linearly with the number of computational nodes. In fact, after some point, increasing the number of computational nodes can lead to slower convergence speed. As mentioned before, this is due to the delay caused by communication among the nodes. This observation gives rise to two important research questions: First, given the maximum allowable number of computation nodes and the  communication overhead of the nodes, what is the optimum choice of the number of cores for solving a  given optimization problem? Answering this question requires computation/communication tradeoff analysis of the proposed optimization approach. Second, can the BSUM framework  be extended and implemented in a (semi-)asynchronous  manner? If this is possible, then the communication overhead  can be reduced significantly since the nodes are not required to wait for each other before updating the variables, making the algorithm {\it lock-free}. For recent efforts on this research direction see \cite{liu14asynchronous}.

\item {\bf Nonlinear coupling constraints}:
As we observe in Section~\ref{sec:extension}, the BSUM framework can also be used in the presence of linear coupling or nonconvex decoupled constraints.  How can the BSUM framework be generalized to problems with nonlinear coupling constraints? More precisely, can the BSUM framework with block-wise update rules be applied to  the optimization problem of the following form?
\begin{align}
\min_x \quad & f_0(x_1,\ldots,x_n) \nonumber\\
\st \quad & f_i(x_1,\ldots,x_n) \leq 0, \;\;\;\forall i=1,2,\ldots,n \nonumber.
\end{align}
Example~\ref{ex:diverge:coupling} shows that the naive extension of the BCD approach fails to find the optimal solution even in the convex setting.  A popular approach to tackle the above problem  is to place the constraints in the objective using Lagrange multipliers and update the multipliers iteratively. However, this approach typically leads to double-loop algorithms and requires subgradient steps in the dual space which is known to be slow.
\end{itemize}

\bibliographystyle{IEEEbib}
\bibliography{ref,biblio}

\end{document}